%



\input amstex
\magnification=\magstep 1
\documentstyle {amsppt}
\topmatter
\title {{\it Remarks on $\aleph_1$-metrizable not metrizable first countable
spaces and CWH} \\
Sh-E9} \endtitle
\rightheadtext{$\aleph_1$-metrizable not metrizable}
\author {Saharon Shelah \thanks {\null\newline I would like to thank Alice
Leonhardt for the beautiful typing. \newline
My thanks also go to Franklin Tall for showing me the problem in 1992, for
the participants of the seminar in Jerusalem in July '92 and Fall '93 for
their remarks, and to Mirna D\v{z}amonja for many corrections. \newline
Publ. No. E9 \newline
Done - \S1,\S2,\S3 June of 1992 - Toronto;\S4 - Oct. 1993. \newline
This version typed 6/24/94} \endthanks} \endauthor
\affil {Institute of Mathematics \\
The Hebrew University \\
Jerusalem, Israel
\medskip
Rutgers University \\
Department of Mathematics \\
New Brunswick, NJ USA} \endaffil
\endtopmatter
\document
\bigskip

\subhead {Abstract} \endsubhead  $CWH,CWN$
stand for collectionwise Hausdorff
and collectionwise normal respectively.
We analyze the statement ``there is a $\lambda-CWH$ not $CWH$ first countable
(Hausdorff topological) space".  We prove the existence of such a space
under various conditions, show
its equivalence to: there is a $\lambda$-CWN not
CWN
first countable space and give an equivalent set theoretic statement; the
nicest version we can obtain is in \S4.  
The author had a flawed proof of the existence of such spaces in ZFC,
for some $\lambda > \aleph_1$, in June of 1992; still we decided that
there is some interest in the correct part and some additions.
\newpage

\head {\S0 Introduction} \endhead
\bigskip

We shall deal mainly with first countable topological spaces. \newline
All spaces will be Hausdorff.
\definition{0.1 Definition}  1) A space $X$ is metrizable if the topology
on $X$ is induced by a metric. \newline
2)  A space $X$ is ($< \lambda$)-metrizable if for each $Y \subseteq X$,
$|Y| < \lambda$, the induced topology on $Y$ is metrizable.  Let
$\mu$-metrizable mean ($< \mu^+$)-metrizable. \newline
3)  A space $X$ is CWH (collectionwise Hausdorff) \underbar{if for every}
subspace $Y$ on which the induced topology is discrete (i.e. every subset
is open) \underbar{there} is a sequence $\langle u_y:y \in Y
\rangle$ of pairwise disjoint open subsets of $X$, such that for every
$y \in Y$ we have $y \in u_y$. \newline
4)  A space $X$ is ($< \lambda$)-CWH if for every $Y \subseteq X$ of
cardinality $< \lambda$, $Y$ (with the induced topology) is $CWH$. \newline
$\mu-CWH$ means ($< \mu^+)-CWH$. \newline
5)  A space is $CWN$ (collectionwise normal) when: if $\langle Y_i:i < \alpha
\rangle$ is a sequence of pairwise disjoint subsets of $X$, and each $Y_i$ is
clopen in $X\restriction ( \dsize \bigcup_{j < \alpha} Y_j )$, then we can
find pairwise disjoint open $\langle {\Cal U}_i:i < \alpha \rangle$ in $X$
such that
$Y_i \subseteq {\Cal U}_i$. \newline
6)  A space is $(< \lambda) -^*$ \, $CWN$ if every subspace with $< \lambda$
points is $CWN$ (we use the $*$ because there may be a bound $\alpha
< \lambda$ such that all relevant subspaces are of size $< \alpha$).
\enddefinition
\bigskip

\noindent
$n -^* CWN$ means ($< \mu^+)-^* CWN$.
\subhead {0.2 Question} \endsubhead  $(ZFC)$  1) Are there
$\aleph_1$-metrizable
not metrizable (first countable Hausdorff topological) spaces? \newline
2)  Are there $\aleph_1-CWH$ not $CWH$ first countable spaces? \newline
We shall also consider analogous questions with $\aleph_1$ replaced by any
$\lambda > \aleph_0$. \newline
\underbar{Note}: $\lambda$-metrizable $\Rightarrow \lambda-CWH$.  Also,
metrizable $\Rightarrow CWN \Rightarrow CWH$.
\bigskip

\subhead{0.3 Observation} \endsubhead  1) Assume $X$ is a space with character
$\chi \le \lambda$ (i.e. every point has a neighborhood basis of cardinality
$\le \chi$). \newline
\noindent
Then:
\medskip
\roster
\item "{(a)}"   $X$ is $\lambda-CWH$ \underbar{iff} for every subspace $Y$ of
cardinality $\le \lambda$ on which the induced topology is discrete there is
a sequence $\langle u_y:y \in Y \rangle$ of pairwise disjoint open subsets
of $X$, $y \in u_y$.
\item "{(b)}"  In (a), for any fixed $\mu \le r$, we can restrict
ourselves (on both sides) to discrete subsets of cardinality $\mu$.
\endroster
\medskip

\noindent
2) If $X$ is $CWN$ then $X$ is $CWH$.
\bigskip

\demo{Proof}  1) The implication $\Leftarrow$ is immediate.  For the
implication $\Rightarrow$ assume that $Y \subseteq X,|Y| \le \lambda$ and
$X \restriction Y$ is
the discrete topology.  Let $\langle {\Cal U}^y_i:i < i^y \le \chi \rangle$ be
a neighborhood basis in $X$ for $y \in Y$; choose for $y^1,y^2 \in Y,
i_1 < i^{y^1},i_2 < i^{y^2}$
a point $z[y^1,y^2,i_1,i_2]$ which is in ${\Cal U}^{y^1}_{i_1} \cap
{\Cal U}^{y^2}_{i_2}$, if this intersection is non-empty.  By the assumption
$X \restriction Y_1$ is $CWH$, where \newline
$Y_1 = Y \cup \{ z[y^1,y^2,i_1,i_2]$:
$y^1 \in Y,y^2 \in Y,i_1 < i^{y^1},i_1 < i^{y^2} \}$.\hfill$\square_{0.3}$
\enddemo\newpage

\head {\S1 Analysis of ``$\aleph_1-CHW$ but not $CHW$"} \endhead
\bigskip

\proclaim{1.1 Lemma}  1)  Assume \newline
\medskip
\roster
\item "{$(*)_\lambda$}"  $cf(\lambda) = \aleph_0 < \lambda$,
$\eta_\alpha \in {}^\omega \lambda$ for $\alpha < \lambda^+$, and for each
$\beta < \lambda^+$, we can find pairwise disjoint end segments for
$\langle \eta_\alpha:\alpha < \beta \rangle$ \newline
(e.g. $\exists h_\beta:\beta \rightarrow \omega$ such that
\endroster

$$
\alpha_1 < \alpha_2 < \beta \wedge k > h_\beta(\alpha_1)\wedge k > h_\beta
(\alpha_2) \Rightarrow \eta_{\alpha_1} \restriction k \ne \eta_{\alpha_2}
\restriction k).
$$
\noindent
\underbar{Then} 1) the space ${}^{\omega >} \lambda \cup \{ \eta_\alpha:\alpha
 < \lambda^+ \}$ with the topology given below is
\medskip
\roster
\item "{$(\alpha)$}"  first countable and Hausdorff
\item "{$(\beta)$}"   $\lambda-CWH$, even $\lambda$-metrizable
\item "{$(\gamma)$}"  \text{not } $\lambda^+ - CHW$.
\endroster
\medskip

\noindent
The topology is the obvious one each $\eta \in {}^{\omega >} \lambda$ is
isolated, and for each $\alpha < \lambda^+$, the neighborhood basis of
$\eta_\alpha$ is $\{ \{ \eta_\alpha \restriction
\ell:k < \ell \le \omega \}:k < \omega \}$. \newline
2)  Moreover, the space is not metrizable but is $\lambda$-metrizable.
\endproclaim
\bigskip

\demo{Proof}  Straightforward. \hfill$\square_{1.1}$
\enddemo
\bigskip

\demo{1.2 Conclusion}  1) If the answer to 0.2(1) or 0.2(2) is ``no",
\underbar{then} $(*)_\lambda$ of 1.1 is not true for any $\lambda$. 
\newline
2)  If $(*)_\lambda$ of 1.1 fails for all $\lambda$, then
\medskip
\roster
\item "{$(*)$}"  $cf(\lambda) = \aleph_0 < \lambda \Rightarrow
pp(\lambda) = \lambda^+$
\endroster
\medskip

\noindent
(by [Sh355,1.5A]). \newline
\noindent
3)  If 2)'s conclusion holds, then for every $\lambda$ singular
we have $pp(\lambda) = \lambda^+$.  (By [Sh371,\,1.10] or [Sh371,\,1.10A(6)]
or [Sh355,\,2.4(1)]), hence for $\theta < \mu$ \newline
$\text{cov}(\mu,\theta^+,\theta^+,2) \le \mu^+$ (by [Sh430, \,1.1]). \newline
4)  If 3)'s conclusion holds then:

$(*)$ $\quad$ if $\lambda$ is singular strong limit then
\medskip
\roster
\item "{(a)}"  $2^\lambda = \lambda^+$
\endroster
hence
\roster
\item "{(b)}"  $\diamondsuit^*_{S_\lambda}$ where $S_\lambda = \{ \delta <
\lambda^+:cf(\delta) \ne cf(\lambda) \}$, and $\diamondsuit^*_S$ means that
\newline
there is a 
$\langle {\Cal P}_\delta:\delta \in S \rangle$, ${\Cal P}_\alpha \subseteq
[\alpha]^\lambda$, $|{\Cal P}_\alpha| = \lambda$ such that
\endroster
$$
(\forall X \subseteq \lambda^+)(\exists \text{ club } C)
\biggl[ \dsize \bigwedge_{\delta \in S \cap C}(X \cap \delta) \in {\Cal P}
_\delta \biggr]
$$

\noindent
(by [Sh108] and see there on earlier work of Gregory). \newline
So clearly $\diamondsuit^*_S \and S_1 \subseteq S \Rightarrow
\diamondsuit^*_{S_1}$.\newline
(5)  Not only $pp(\lambda) > \lambda^+$ and $\lambda > \aleph_0 =
cf(\lambda)$ implies $(*)_\lambda$ (from 1.11); but assume we have
$\langle \lambda_n:n < cf(\lambda) \rangle$, $\dsize \sum_n \lambda_n =
\lambda$, $\lambda_n = cf(\lambda_n)$,
tcf$(\Pi \lambda_n / J^{bd}_\omega) = \lambda^+$ exemplified by
$\bar f = \langle f_\alpha:\alpha < \lambda^+ \rangle$ such that
\medskip
\roster
\item "{$\oplus$}"  if $\aleph_0 < cf(\delta) = \kappa < \lambda$, then
there is a closed unbounded $A \subseteq \delta$ and $n_\alpha <
cf(\lambda)$ for $\alpha \in A$ such that $n_\alpha,n_\beta < n <
cf(\lambda) \Rightarrow f_\alpha(n) < f_\beta(n)$.
\endroster
\medskip

\noindent
\underbar{Then} using $\oplus$ rather than $(*)_\lambda$, in 1.1 we get a
$\kappa^+$-CWH, $\kappa^+$-metrizable first
countable space (see [Sh355,\S6]). \hfill$\square_{1.2}$
\enddemo
\bigskip

\demo {1.3 Construction}  Assume $\lambda = \beth_{\omega_1}$ (or just
$\lambda$ is a strong limit, $cf(\lambda) \ne \aleph_0)$,
$2^\lambda = \lambda^+$ and
$S$ is a stationary subset of $\lambda^+$, \newline
$S \subseteq \{ \delta < \lambda^+:
cf(\delta) = \aleph_0 \text{ and } \omega^2 \text{ divides } \delta \}$
(the existence of an $S$ like that such that
$\diamondsuit_S$ suffices). \newline
We shall build a space with the set of points
$\{ x_\alpha,y_\alpha:\alpha < \lambda^+ \}$.  Each $x_\alpha$ will be
isolated in $X$ and each $y_\beta$ will have a countable neighborhood basis
in $X$.  We shall have $\{ u_{\alpha,n}:n < \omega \}$ as a neighborhood base
of $y_\alpha$ with $u_{\alpha,n}$ decreasing in $n$ and
$u_{\alpha,n} = \{ y_\alpha \} \cup \{ x_\beta:f_\alpha(\beta) > n \}$
where $f_\alpha(\beta) \in \omega$. \newline
Note that each $Y_\alpha$ is isolated in the space restricted to
$\{ Y_\alpha:\alpha < \lambda^+ \}$. \newline
The only thing left is to define is $f$. \newline
We set $f_\alpha(\beta) = 0$ except in some specified cases.  For the space
to be Hausdorff it is enough to have: \newline

for $\alpha < \beta$ there is an $m = m(\alpha,\beta) < \omega$ such that
\newline

$\neg(\exists \gamma)[f_\alpha(\gamma) \ge m \and
f_\beta(\gamma) \ge m]$.  We shall make a stronger condition:
\medskip
\roster
\item "{$(*)$}"  $\alpha < \beta \Rightarrow (\exists^{\le 1} \gamma)
[f_\alpha(\gamma) \ge 1 \and f_\beta(\gamma) \ge 1]$.
\endroster
\medskip

\noindent
Remember that $\diamondsuit_S$ holds as $2^\lambda = \lambda^+$ and
$cf(\lambda) > \aleph_0$.  So there is a $\langle g_\alpha:\alpha \in S
\rangle$ \,\, $g_\alpha:\alpha \rightarrow \omega$ such that

$$
(\forall g \in {}^{\lambda^+}\omega)(\exists^{\text{stat}}\alpha \in S)
(g_\alpha = g\restriction \alpha).
$$

\noindent
Now, if the space is CWH then there is a $g:\lambda^+ \rightarrow \omega$
such that
$\langle u_{\alpha,g(\alpha)}:\alpha < \lambda^+ \rangle$ are pairwise
disjoint.

We define by induction on $\alpha$ a limit $< \lambda^+$, $f_i(j)$ for
$i,j < \alpha$.  Call the sequence $\langle f_i:i < \alpha \rangle$ in
${}^\alpha \omega$ \, $f^\alpha$, so if $\alpha < \beta$, then $f^\alpha$
is an initial segment of $f^\beta$.  Usually we just give value zero to
$f^i(j)$. \newline
\noindent
If $\alpha \in S$, and $g_\alpha$ looks as a candidate for $g$, i.e.
$\langle u^\alpha_{i,g_\alpha(i)}:i < \alpha \rangle$ are pairwise disjoint,
where
$u^\alpha_{i,k} =: \{ \beta < \alpha:f_i(\beta) > k \}$, and if
for some $m = m_\alpha$, \newline
otp$(\{ \beta < \alpha:g_\alpha(\beta) = m \}) = \alpha$, then choose
\medskip
\roster
\item "{(a)}"  $\beta^\alpha_n < \beta^\alpha_{n+1} \dots < \alpha =
\dsize \bigcup_{n} \beta^\alpha_n$
\item "{(b)}"  $g_\alpha(\beta^\alpha_n) = m$
\endroster
\medskip

\noindent
and define $f^{\alpha + \omega}$ (extending $f^\alpha$) by

$$
f^{\alpha + \omega}_\alpha(\alpha + n) = n
$$

\noindent
and

$$
f^{\alpha + \omega}_{\beta^\alpha_n}(\alpha + n) = m + 1
$$

\noindent
(other values of $f^{\alpha + \omega}$ are zero).  If $g_\alpha$ fails the
conditions above, choose $m_\alpha = 0$, $\beta^\alpha_n$ satisfying
conditions (a) above and extend $f^\alpha$ as just described. \newline
So we cannot extend $g_\alpha$ to $n^+:\text{if } g(\alpha) = k$ we
get

$$
\in (u_{\alpha,k+1} \cap u_{\beta^\alpha_{k+18}},g_\alpha(\beta^\alpha
_{k+18})).
$$

\noindent
So the space is not CWH (hence not metrizable).  For simplicity, we can
request that
$\beta^\alpha_n \notin \dsize \bigcup_{\gamma \in S} [\gamma,\gamma +
 \omega)$.  Suppose the space is not $\aleph_1-CWH$.  So for some
${\Cal U} \in [\lambda^+]^{\aleph_1}$

$$
X\restriction \{ x_\alpha,y_\alpha:\alpha \in {\Cal U} \}
$$

\noindent  is not $CWH$. \newline
So without loss of generality if

$$
\alpha \in S \cap {\Cal U}
$$
\noindent
then

$$
\gather
\alpha + n \in {\Cal U} \text{ and} \\
\beta^\alpha_n \in {\Cal U}.
\endgather
$$
\medskip
\noindent  So
\roster
\item "{$\otimes$}"  for every $g:{\Cal U} \rightarrow \omega$ (candidate to
give the separation), \newline
\underbar{we get}: for some $\alpha \in S \cap {\Cal U}$, $(\exists^\infty n)
\,g(\beta^\alpha_n) \le m_\alpha$.
\endroster

\noindent
This is a contradiction. \hfill$\square_{1.3}$
\enddemo
\bigskip

\subhead {1.4 Comments} \endsubhead
\bigskip

\noindent
(1)  The space constructed in 1.3 does not have neighborhood bases
consisting of countable sets,  so
is not excluded by the earlier consistency results from [JShS320]. \newline
(2)  But $\Vdash_{\text{Levy}(\aleph_1,\lambda^+)}$ ``$X$ is not $\aleph_1$-
CWH" may fail unless we put more restrictions on the $\beta^\alpha_n$.
See (3). \newline
(3)  \underbar{If} we build $X$ as above, let $P = \text{Levy}(\aleph_1,
\lambda^+)$ and we build a $P$-name ${\underset\sim\to g}$ such that

$$
\Vdash_P ``{\underset\tilde {}\to g}:\lambda \rightarrow \omega
 \text{ witnesses that } X \text{ is } CWH",
$$

\noindent
then $X$ is $\aleph_1$-CWH. \newline
[Why?  given a $Y \in [\lambda^+]^{\aleph_1}$, we can find
$\langle p_i:i < \omega_1 \rangle$ increasing in $P$ such that \newline
$\dsize \bigwedge_{\alpha \in Y}
\dsize \bigvee_{i} p_i \Vdash {\underset\tilde {}\to g}(\alpha) = $
 something]. \newline
\bigskip

\definition{1.5 Definition}  We say that the space $X$ is $\lambda-WCWH$
if for any discrete set of $\lambda$ points, some subset of
cardinality $\lambda$ can be separated by disjoint open sets.
\enddefinition
\bigskip

\remark{1.5A Remark}  By a theorem of Foreman and Laver for first countable
spaces we have the consistency of:
$\aleph_1-WCWH \Rightarrow \aleph_2-WCWH$.

On the other hand, e.g. namely in [FoLa], starting with a huge embedding
$j:V \rightarrow M$ with critical point $\kappa$ and $j(\kappa) = \lambda$,
the following is obtained:

There is a forcing notion $\Bbb P * \Bbb R$ such that $\Bbb P$ is
$\kappa$-c.c.,$| \Bbb P| = \kappa,V[G_{\Bbb P}] \models \kappa = \omega_1$,
$\Bbb R \in V[G_{\Bbb P}]$ is $\lambda$-c.c., of cardinality $\lambda$ and
$(< \kappa)$-closed and $V[G_{\Bbb P * \Bbb R}] \models \lambda =
\omega_2$.  In addition, there is a regular embedding $h:(\Bbb P * \Bbb R)
\rightarrow j,\Bbb P$ with $h(p) - p$ for all $p \in \Bbb P$ and the master
condition property holds for $h,j\Bbb P * \Bbb R$.  Finally, if $G$ is
$(\Bbb P * \Bbb R)$-generic, then in $V[G],j \Bbb P / h''(G)$ is $\kappa$-
centered.

The consistency of $\aleph_1-WCWH \rightarrow \aleph_2-WCWH$ for first
countable spaces clearly follows from the above result of [FoLa].  For the
convenience of the reader we include the following easy Claim 1.5B which
shows this implication.
\endremark
\bigskip

\proclaim {1.5B Claim}  Suppose $X$ is a first countable topological space
and $|X| = \kappa^+$, while $Y_0 \subseteq X$ is a discrete subspace of $X$,
with $|Y_0| = \kappa^+$.  If $\Bbb P$ is a $\kappa^+$-c.c., even
$\kappa$-centered forcing notion such that

$$
\Vdash_P ``\text{There is a } Y \subseteq Y_0 \text{ with }
|Y| = |X| \text{ and } Y \text{ is separated in } X",
$$

\noindent
\underbar{then} \newline
in $V$, there is a $Y \subseteq Y_0,|Y| = |Y_0|$ and $Y$ is separated in
$X$.
\endproclaim
\bigskip

\demo{Proof}  Without loss of generality, the set of points of $Y_0$ in
$V^P$ is $\kappa^+$, and we denote $\lambda = \kappa^+$.  We may fix a set
$\{ {\underset\tilde {}\to x_\gamma}:\gamma < \lambda \}$ of
$\Bbb P$-names such that

$$
\Vdash_P ``\{ {\underset\tilde {}\to x_\gamma}:\gamma < \lambda \}
\text{ is separated}".
$$
We can also assume that there are no repetitions among the
${\underset\tilde {}\to x_\gamma}$, and that
${\underset\tilde {}\to x_\gamma} \ge \gamma$. \newline
Suppose that in $V$, the neighborhood bases for points in $Y_0$ are given
by

$$
\left< \langle u^n_y:n < \omega \rangle:y \in Y_0 \right>.
$$
So, without loss of generality $\{ {u^{\underset\tilde {}\to n(\gamma)}
_{\underset\tilde {}\to x_\gamma}}:\gamma < \lambda \}$ are pairwise
disjoint, in $V^{\Bbb P}$.

Now, let $\Bbb P = \dsize \bigcup_{i < \kappa} \Bbb P_i$ where each
$\Bbb P_i$ is directed. \newline
For each $\alpha < \lambda$, there is a forcing value to
${\underset\tilde {}\to x_\alpha}$, say $\beta_\alpha$.  So, there is an
$i(*) < \kappa$ such that $A = \{ \alpha:\beta_\alpha \in P_{i(*)} \}$ is
unbounded in $\lambda$. \newline
Therefore, $\{ \beta_\alpha:\alpha \in A \}$ is separated by

$$
\{ u^{n(\alpha)}_{\beta_\alpha}:\alpha \in A \}.
$$

\noindent
(So, having that any two members of $\Bbb P_i$ are compatible, or that out
of any $\lambda$ elements of $\Bbb P$ there are $\lambda$ pairwise
compatible, i.e. $\Bbb P$ is $\lambda$-Knaster, suffices).
\hfill$\square_{1.5.B}$
\enddemo
\bigskip

\proclaim{1.6 Claim}  There is a first countable Hausdorff space $X$ which is
$(2^{\aleph_1})^+-WCWH$ but is not $WCWH$.
\endproclaim
\bigskip

\demo{Proof}  Let $\lambda = \dsize \sum_{n < \omega} \lambda_n$,
$\lambda^{\aleph_0}
_n < \lambda_{n+1}$.  Let $\langle \eta_\alpha:\alpha < \lambda^+ \rangle$,
$\eta_\alpha \in {}^\omega \lambda$, $\alpha < \beta$ and $\eta_\alpha
<_{J^{bd}_\omega} \eta_\beta$. \newline
Topology: as in 1.1. \newline
\underbar{Proof of not $\lambda^+-WCWH$}: if ${\Cal U} \in
[\lambda^+]^{\lambda^+}$,
$\langle \eta_\alpha:\alpha \in {\Cal U} \rangle$ cannot be separated \newline
as $|\{ \eta_\alpha \restriction \ell:\ell < \omega,\alpha \in {\Cal U} \}|
\le \lambda$.

If ${\Cal U} \in [\lambda^+]^{(2^{\aleph_0})^+}$, without loss of generality
otp$({\Cal U}) = (2^{\aleph_0})^+$; set \newline
${\Cal U} = \{ \alpha_\zeta:\zeta < (2^{\aleph_0})^+ \}$.  Now for some
$Y \in [(2^{\aleph_0})^+]^{(2^{\aleph_0})^+}$ and $n$, \newline
$\langle \eta_{\alpha_\zeta} \restriction [n,\omega):\zeta \in Y
\rangle$ is strictly increasing (not just modulo $J^{bd}_\omega$ but
in every coordinate (see [Sh355,\S6], [Sh400,\S5], [Sh430,\S6]).
\hfill$\square_{1.6}$
\enddemo
\bigskip

\remark{1.7 Remark}  We can prove other Claims like 1.6 (see the
references above).
\endremark
\newpage

\head {\S2 On not CWH, $\aleph_1$-CWH Spaces} \endhead
\bigskip

\definition{2.1 Definition}  For an ordinal $\gamma$ let us define
\medskip
\roster
\item "{$(*)^1_\gamma$}"  there is a $S \subseteq \{ \delta < \gamma:
cf(\delta) = \aleph_0 \}$ and, for $\delta \in S$, a sequence
$\langle \beta^\delta_n:n < \omega
\rangle$ strictly increasing with limit $\delta$, and a $m_\delta < \omega$,
\newline
such that $(\forall g \in {}^\gamma\omega)(\exists \delta \in S)
(\exists^\infty n)[g(\beta^\delta_n) \le m_\delta]$.
\endroster
\enddefinition
\bigskip

\proclaim{2.2 Claim}  (1)  If the answer to 0.2 is no (or much less), then
for some $\gamma < \omega_2$, $(*)^1_\gamma$ holds. \newline
(2)  If $MA + \neg CH$, then
$\gamma < 2^{\aleph_0} \Rightarrow \neg(*)^1_\gamma$.
\newline
(3) Without loss of generality, in $(*)^1_\gamma$, each $\beta^\delta_n$ is a
successor ordinal.  
\endproclaim 
\bigskip

\demo{Proof}  1) By the proof of 1.3 and 1.2. \newline
(2)  Check.  Use the natural forcing $\{ p:p$ is a finite function from
$\gamma$ to $\omega \}$ with $p \le g$ iff $p \subseteq g \and
(\forall \delta)(\delta \in S \cap \text{ Dom}(p) \rightarrow (\forall n)
[\beta_n \in \text{ Dom}(g) \backslash \text{Dom}(p) \rightarrow g(\beta_n)
> n_\delta]$). \newline
(3)  Check.
\enddemo
\bigskip

\demo {2.2A Conclusion}  If $MA + \neg CH$ \underbar{then} the answer to 0.2
is yes.  In fact, there is an $\aleph_1$-metrizable (hence $\aleph_1$-CWH) not
CWH (hence not metrizable) first countable space.
\enddemo
\bigskip

\demo{Proof}  By 2.2(1) and 2.2(2).
\enddemo
\bigskip

\proclaim{2.3 Claim}  If $(*)^1_\gamma$ for some
$\gamma < \omega_2$, then $(*)^1_{\omega_1}$.
\endproclaim
\bigskip

\demo{Proof}  Choose $\gamma^* < \omega_2$ minimal such that
$(*)^1_{\gamma^*}$.
Clearly $\gamma^* \ge \omega_1$.

If $\gamma^* = \omega_1$ we are done.  So assume $\gamma^* > \omega_1$,
and we shall get a contradiction.  We fix an $S \subseteq \gamma$ and $mJ$,
$\langle \beta^J_n:n < \omega \rangle$ for $J \in S_1$, which exemplify
$(*)^1_{\gamma^*}$.  Note that for every $\gamma < \gamma^*$ there is a
$g_\gamma \in {}^\gamma\omega$ such that:
\medskip
\roster
\item "{$\otimes$}"  if $\delta \in S \cap \gamma$ \underbar{then}
$\{ n:g_\gamma(\beta^\delta_n) \le m_\delta \}$ is finite.
\endroster\enddemo
\bigskip

\subhead {Case 1} \endsubhead  $\gamma^* = \gamma +1$, $\gamma \notin S$.
\newline
Extend $g_\gamma$ by $\{ \langle \gamma,0 \rangle \}$.
\bigskip

\subhead {Case 2} \endsubhead  $\gamma^* = \gamma +1$, $\gamma \in S$:\newline
define $g \in {}^{\gamma^*}\omega$:

$$
\alignat2
\text{\underbar{if} } &\beta \in \gamma,\beta \notin \{
\beta^\gamma_n:n < \omega \} \qquad &&\text{\underbar{then} } g(\beta) =
g_\gamma(\beta) \\
\text{\underbar{if} } &\beta = \gamma \qquad &&\text{\underbar{then} }
g(\beta) = 0 \\
\text{\underbar{if} } &\beta = \beta^\gamma_n \qquad &&\text{\underbar{then} }
g(\beta) = \text{ Max} \{ g_\gamma(\beta),n+8,m_\gamma + 8 \}.
\endalignat
$$

\noindent
So $g$ gives a contradiction.
\bigskip

\subhead {Case 3} \endsubhead  $cf(\gamma^*) = \aleph_0$.

Let $\gamma^* = \dsize \bigcup_{n < \omega} \gamma_n$, $\gamma_0 = 0$,
$\gamma_n < \gamma_{n+1}$, and each $\gamma_{n+1}$ is a successor of a
successor ordinal.

Let $g = \cup \{ g_{\gamma_{n+1}} \restriction [\gamma_n,\gamma_{n+1}):n <
\omega \}$ - it gives a contradiction.
\bigskip

\subhead {Case 4} \endsubhead  $cf(\gamma^*) = \omega_1$.

Let $\langle \gamma_i:i < \omega_1 \rangle$ be increasing continuous with
limit $\gamma^*$, $\gamma_0 = 0$, $\gamma_{i+1}$ a successor of a successor
ordinal. \newline
Let $S' =: \{ \gamma_i:\gamma_i \in S \text{ (so } i \text{ is a limit
ordinal)} \}$.
\bigskip

\subhead {Subcase A} \endsubhead  $\gamma^*$, $\langle <\beta^\gamma_n:
n < \omega>:\gamma \in S' \rangle$, $\langle m_\gamma:\gamma \in S'
\rangle$ does not exemplify $(*)^1_{\gamma^*}$.  So some $g^* \in
{}^{\gamma^*}\omega$ shows this.  Define $g$ by:

$$
\text{\underbar{if} } \beta \in [\gamma_i,\gamma_{i+1}) \,
\text{\underbar{then }}
 g(\beta) = \text{ Max}
\{g_{\gamma_{i+1}}(\beta),g^*(\beta)\}
$$

\noindent
So $g$ gives a contradiction.
\bigskip

\subhead {Subcase B} \endsubhead  $\langle <\beta^\gamma_n:n < \omega>:
\gamma \in S' \rangle$, $\langle m_\gamma:\gamma \in S' \rangle$
exemplifies $(*)^1_{\gamma^*}$. \newline
Let $S^* = \{ i < \omega_1:i \text{ a limit},\gamma_i \in S \}$
(necessarily stationary). \newline
Let $\gamma^* = \dsize \bigcup_{i<\omega_1} a_i$, $a_i$ countable increasing
continuous, \newline
\underbar{such that} $a_0 = \emptyset$, \newline
$a_i \cap \{ \gamma_j:j < \omega_1 \} = \{ \gamma_j:j < i \}$, 
$a_i \subseteq \gamma_i$ and $\gamma_j \in a_i \wedge j \in S^* \Rightarrow
\dsize \bigwedge_{n} \beta^{\gamma_j}_n \in a_i$. \newline
For $i \in S^*$ let $u_i = \{ n < \omega:\beta
^{\gamma_i}_n \in a_i \}$. \newline
Note
\medskip
\roster
\item "{$\otimes$}"  if $i \in S^*$ and $j < i$, then $\{ n \in u_i:\beta
^{\gamma_i}_n \in a_j \}$ is finite, as it is included in $\{ n < \omega:
\beta^{\gamma_i}_n < \gamma_j \}$.  [Why?  Remember $a_j \subseteq \gamma_j$].
\endroster
\medskip

\noindent
Let $S^{**} = \{ i \in S^*:u_i \text{ is infinite and } i \text{ is a limit
ordinal}\}$.  So we already know
\medskip
\roster
\item "{$\oplus$}"  for every $g \in {}^{\gamma^*}\omega$, for some
$i \in S^*$, for infinitely many \newline
$n < \omega$, $g(\beta^{\gamma_i}_n) \le m_{\gamma_i}$.
\endroster
\medskip

\noindent
We claim
\medskip
\roster
\item "{$\oplus^+$}"  for every $g \in {}^{\gamma^*}\omega$ for some
$i \in S^{**}$, for infinitely many $n \in u_i$ we have $g(\beta^{\gamma_i}
_n) \le m_{\gamma_i}$.
\endroster
\medskip

\noindent
Otherwise, for some $g^* \in {}^{\gamma^*}\omega$ this fails and we
define $g$:
\roster
\item "{{}}"  \underbar{let} $\beta \in a_{i+1}\backslash a_i$ (there is one
and only one such $i$), \newline
\underbar{then} $g(\beta) = \text{ Max}\{ g^*(\beta),m_{\gamma_i} + 8,
m_{\gamma_i+1} + 8 \}$
\endroster
\medskip

\noindent
As $g$ gives a contradiction to $\oplus$, clearly $\oplus^+$ holds.

Now let $h$ be a one to one function from $\omega_1$ onto $\gamma^*$ such that
for $i$ limit, $h$ maps $\{j:j < i \}$ onto $a_i$. \newline
Let for $i \in S^{**}$, $\{ j^i_n:n < \omega \}$ enumerate $\{j < i:h(j) \in
\{ \beta^{\gamma_i}_n:n \in u_i \} \}$, and $m^*_i = m_{\gamma_i}$ for
$i \in S^{**}$. \newline
Now $\langle <j^i_n:n < \omega>:i \in S^{**} \rangle$, $\langle m^*_i:i \in
S^{**} \rangle$ exemplifies that $\gamma^*$ could have
been chosen to be $= \omega_1$, as required.\hfill$\square_{2.3}$
\medskip

\noindent
We define the combinatorial property we actually use 
\definition{2.4 Definition}  1) $INCWH(\lambda) = INCWH^1(\lambda)$ means:
\medskip
\roster
\item "{{}}"  $\lambda$ is regular $> \aleph_0$ and for some stationary
$S \subseteq \{ \delta < \lambda:cf(\delta) = \aleph_0 \}$ we have
$\langle m_\delta < \beta^\delta_n:n < \omega>:
\delta \in S \rangle$ such that:
\newline
$m_\delta < \omega$, $\beta^\delta_n < \beta^\delta_{n+1} < \delta =
\dsize \bigcup_{n < \omega} \beta^\delta_n,\beta^\delta_n$ is a
successor \underbar{and}:
{\roster
\itemitem{ $(a)$ }  for every $g \in {}^\lambda \omega$, for some $\delta \in
S$, for infinitely many $n$, $g(\beta^\delta_n) \le m_\delta$
\itemitem{ $(b)_\lambda$ }  for every $U \subseteq \lambda$,
$|U| < \lambda$, for some $g \in {}^U \omega$, for every $\delta \in S
\cap U$, for every $n < \omega$ large enough, 
$g(\beta^\delta_n) > m_\delta$.
\endroster}
\endroster

\noindent
2) We can replace $m_\delta$ by $\langle m^\delta_n:n < \omega \rangle$,
requesting $q(\beta^\delta_n) \le m^n_\delta$ in (a) and
$q(\beta^\delta_n) > m^n$ in (b)$_\lambda$.  In this way we obtain a
weaker property, which we call $INCWH^2(\lambda)$. \newline
For other versions of the principle, as well as the connections between the
various versions, see \S3.
\enddefinition
\bigskip

\demo{2.4A Discussion}  1) If $INCWH(\lambda)$, then there is a space (as in
1.3) which is Hausdorff first countable with $\lambda$ points, not
metrizable, not even CWH, but every subspace of smaller cardinality
is metrizable. \newline
2)  So if we prove $(\exists \lambda > \aleph_1)INCWH(\lambda)$ we have
solved the original problem 0.2. \newline
3)  $(b)_\kappa$ means that we require $|U| < \kappa$.  Note that
(b)$_{\aleph_1}$ holds trivially and that $n \le \kappa \and (b)_\kappa
\Rightarrow (b)_n$.
\enddemo
\bigskip

\noindent
More formally
\proclaim{2.5 Claim}  If $INCWH(\lambda)$ \underbar{then} $SINCWH(\lambda)$
(even exemplified by a $(< \lambda)$-metrizable space) where:
\endproclaim
\bigskip

\definition{2.6 Definition}  $SINCWH(\lambda)$ means that
there is a first countable
$T_2$-space $X$ with $\lambda$ points which is ($< \lambda$)-CWH (i.e. for
every discrete subset of cardinality $< \lambda$ we can choose pairwise
disjoint open neighborhoods) but not $\lambda$-CWH.
\enddefinition
\bigskip

\demo{Proof of 2.5}  The points of $X$ are $y_\alpha$ \,\,
$(\alpha < \lambda)$ and $x_{\alpha,\beta}$ \,\,
$(\beta < \alpha < \lambda)$ with $x_{\alpha,\beta}$
isolated and $y_\alpha$ which have neighborhood bases
$\langle u_{\alpha,n}:n < \omega \rangle$:

$$
\text{\underbar{if}} \quad
 \alpha \in S \qquad u_{\alpha,n} = \{y_\alpha\} \cup
\{x_{\alpha,\beta}:\text{for some } k > n,\beta = \beta^\alpha_k \}
$$

$$
\text{\underbar{if}} \,\, \alpha \notin S \qquad u_{\alpha,n} =
 \{y_\alpha\} \cup
\{x_{\delta,\alpha}:\alpha < \delta \in S, \text{ for some } k,\alpha
 = \beta^\delta_k, \text{ and } n \le m^\delta \}.
$$

\noindent
Here, $S$ is a fixed stationary $\subseteq \{ \delta < n:cf(\delta) =
\aleph_0 \}$ which exemplifies $INCWH(\lambda)$, together with
$\left< m_\delta,\langle \beta^\delta_n:n < \omega \rangle:\delta \in S
\right>$. 
\enddemo
\bigskip

\noindent
\underbar{Checking of ``$X$ not $CWH$"} \newline

Let $Y = \{ Y_\alpha:\alpha < \lambda \}$. \newline
Note that $X \restriction Y$ is a discrete subspace of $X$.  Let
$\{ u_{\alpha,n}:n < \omega \}$ be the neighborhood basis of $y_\alpha$
there is
$\langle u_{\alpha,g(\alpha)}:\alpha < \lambda \rangle$, a sequence of
pairwise disjoint sets, for some $g \in {}^\lambda \omega$.  As
$u_{\alpha,g(\alpha)} \cap u_{\beta,g(\beta)} = \emptyset$ for
$\alpha \ne \beta(< \lambda)$ clearly for
$\alpha \in S,\beta = \beta^\delta_n$ we get $n > g(\alpha) \Rightarrow
g(\beta) > m_\alpha$ (since otherwise $x_{\alpha,\beta} \in
u_{\alpha,g(\alpha)} \cap u_{\beta,g(\beta)}$. \newline
So $g$ contradicts (a) of $INCWH(\lambda)$.
\bigskip

\underbar{Checking of ``$X$ is $(< \lambda)-CWH$"} \newline

Let $Z \subseteq X,|Z| < \lambda$ and $X \restriction Z$ is discrete.
Let \newline
$Z_0 = \{x_{\alpha,\beta}:\beta < \alpha < \lambda \} \cap Z,Z_1 = \{y_\alpha:
\alpha \in \lambda \backslash S \} \cap Z$,
$Z_2 = \{ y_\alpha:\alpha \in S \} \cap Z$, so
$\langle Z_1,Z_2,Z_3 \rangle$ is a
partition of $Z$.  Let ${\Cal U} = \{ \alpha \in S:y_\alpha \in Z_2 \}$,
so $| {\Cal U} | < \lambda,{\Cal U} \subseteq \lambda$ hence by the
assumption,
there is a $g_0 \in {}^\lambda \omega$ as in $(b)_\lambda$. \newline
We define $u_z$, a neighborhood of $z$ for $z \in Z$:
\medskip
\roster
\item "{{}}"  if $z \in Z_0,u_z = \{ x_{\alpha,\beta} \}$
\item "{{}}"  if $z = y_\alpha \in Z_1,u_z = u_{\alpha,n(\alpha)}$
 where \newline

$\qquad \qquad n(\alpha) = \text{ Min}\{n:n \ge g(\alpha) + 8 \text{ and }
u_{\alpha,n} \cap Z_0 = \emptyset \}$
\item "{{}}"  if $z = y_\delta \in Z_2,u_z = u_{\delta,n(\delta)}$ where
\newline

$\qquad \qquad n(\delta) =
 \text{ Min}\{n:n \ge m_\delta + 8 \text{ and } u_{\delta,n} \cap Z_0
 = \emptyset \}$.
\endroster

\noindent
Now check. \hfill$\square_{2.5}$
\bigskip

\proclaim{2.7 Claim}  Assume $\lambda$, $\left< m_\delta,\langle
\beta^\delta_n:n < \omega \rangle:\delta \in S \right>$,
are as in 2.4 but we require $\lambda$ just to be an ordinal, and weaken
(b)$_\lambda$ to
\medskip
\roster
\item "{$(b)_\kappa$}"  for every ${\Cal U} \subseteq \lambda$, $|{\Cal U}| <
\kappa$, for some $g \in {}^U \omega$ \newline
for every $\delta \in S \cap U$, for every $n$ large enough $g(\beta^\delta_n)
> m_\delta$.
\endroster
\medskip

\noindent
\underbar{Then} for some regular $\mu$, $\kappa \le \mu \le \lambda$ we have
$INCWH(\mu)$.
\endproclaim
\bigskip

\demo{Proof}  If we allow $\mu$ in the definition of $INCWH(\mu)$ to be an
ordinal: straightforward (and suffices for
our main interest).  Namely, we choose a ${\Cal U}$ such that
\medskip
\roster
\item "{$(\alpha)$}"  ${\Cal U} \subseteq \lambda$,
\item "{$(\beta)$}"  there is no $g \in {}^\lambda \omega$ such that for every
$\delta \in S \cap {\Cal U}$ for every $n$ large enough $g(\beta^\delta_n) >
 m_\delta$,
\item "{$(\gamma)$}"  under $(\alpha) + (\beta)$ the order type of ${\Cal U}$
is minimal.
\endroster
\medskip

\noindent
Clearly otp$({\Cal U}) \le \lambda$ and otp$({\Cal U}) \ge \kappa$.  By the 
same proof of 2.3, otp$(U)$ is a regular cardinal, we call it $\mu$ and with
the $a_i$'s as in the proof of 2.3, we get $INCWH(\mu)$. \hfill$\square_{2.7}$
\enddemo
\bigskip

\demo{2.8 Conclusion}  If $\lambda = cf(\lambda) > \aleph_0$, $\diamondsuit
_{\{\delta < \lambda:cf(\delta) = \aleph_0\}}$ then for some regular
uncountable $\lambda' \le \lambda$, $INCWH(\lambda')$.  This follows by the
proof of 1.3 and (b)$_{\aleph_1}$ and 2.7.
\enddemo
\bigskip

\demo{2.9 Observation}  If $S_1 \subseteq S_2 \subseteq \{ \delta <
\lambda:cf(\delta) = \aleph_0 \}$, $\left< m_\delta,\langle
\beta^\delta_n:n < \omega \rangle:\delta \in S_1 \right>$ witness
$INCWH(\lambda)$, \underbar{then} we can find a $\left< m'_\delta,
\langle \gamma^\delta_n:n < \omega \rangle:\delta \in S_2 \right>$
witnessing $INCWH(\lambda)$.
\enddemo
\bigskip

\remark{2.10 Remark}  1) We can replace in our discussion
$\aleph_0$ by $\theta$.  Toward this we define a family of spaces.
\endremark
\bigskip

\definition{2.11 Definition}  $X \in {\Cal T}^\ell_\theta$ \underbar{if}
$X$ is a Hausdorff space with each point $x$ having a neighborhood basis
$\{u_{x,\alpha}:\alpha < \alpha^*\}$
such that:
\roster
\item {(a)}  $\ell = 0$ and $\alpha^* \le \theta$
\item {(b)}  $\ell = 1, \alpha^* \le \theta$ and $\langle u_{x,\alpha}:
\alpha < \alpha^* \rangle$ is decreasing.
\item {(c)}  $\ell = 2, \alpha^* = \theta$, and $\langle u_{x,\alpha}:\alpha
< \alpha^* \rangle$ is decreasing.
\endroster
\enddefinition
\bigskip

\definition{2.12 Definition}  We define also the principles \newline
$INCWH(\lambda,\theta) = INCWH^1(\lambda,\theta)$ and
$INCWH^2(\lambda,\theta)$ as in 2.4.
\enddefinition
\bigskip

\proclaim{2.13 Claim}
\roster
\item "{$(\alpha)$}"  if $\lambda > cf(\lambda) = \theta$, $pp(\lambda) >
\lambda^+$ (or the parallel of 1.2(5)), then
{\roster
\itemitem{ $\otimes$ }  there is an $X \in {\Cal T}^2_\theta$,
$|X| = \lambda^+$,
$X$ is $\lambda-CWH$, $X$ has a discrete subspace of size $\lambda^+$,
but for some $X' \subseteq X$, $|X'| = \lambda$,
$cl(X') = X$ (so $|cl(X')| > \lambda)$ (this is a strong form of $X$
is not $\lambda^+-CWH$).
\endroster}
\item "{$(\beta)$}"  if $\lambda > cf(\lambda) = \theta$, $\lambda$ is a
strong limit and $2^\lambda = \lambda^+$, then:
$INCWH(\lambda',\theta)$ for some
$\lambda' = cf(\lambda') \in [\theta^+,\lambda^+]$.
\endroster
\endproclaim
\bigskip

\demo{Proof}  Similar to the above. \hfill$\square_{2.13}$
\enddemo
\newpage

\head {\S3 Variants of Freeness} \endhead
\bigskip

\definition{3.1 Definition}  1) $INCwh(\lambda) = INCwh^1(\lambda)$
is defined as in 2.4 except
that $\langle \beta^\delta_n:n < \omega \rangle$ is not required to be
increasing with limit $\delta$, just $[n \ne m \Rightarrow \beta^\delta_n
\ne \beta^\delta_m]$. \newline
2)  $INCwh^2(\lambda)$ is defined as in (1) but we use $\langle m^\delta_n:
n \in \omega \rangle$ rather than a single $m_\delta$.
\enddefinition
\bigskip

\proclaim{3.2 Claim}  0) $INCWH^\ell(\lambda) \Rightarrow INCwh^\ell
(\lambda)$, $INCWH^1(\lambda) \Rightarrow INCWH^2(\lambda)$,\newline
$INCwh^1(\lambda) \Rightarrow INCwh^2(\lambda)$. \newline
1)  $INCwh^2({\frak b})$ (where \newline
${\frak b} = \text{ Min} \{|F|:f \subseteq {}^\omega \omega
 \text{ and for no } g \in {}^\omega \omega \text{ for every } f \in F,
f <^* g \}$. \newline
2)  Assume $\lambda \le 2^{\aleph_0}$ and for $\alpha < r$,
$f_\alpha$ is a partial function from $\omega$ to $\omega$,
$\text{Dom}(f_\alpha)$ is infinite and $U \subseteq \lambda \and |u|
< \lambda \Rightarrow (\exists f \in {}^\omega \omega) \dsize \bigwedge
_{\alpha \in {\Cal U}} f_\alpha \le^* f$ but for no $f \in {}^\omega \omega$,
$\dsize \bigwedge_{\alpha < \lambda} f_\alpha <^* f$, \newline
\underbar{then} $INCwh^2(\lambda)$. \newline
3) It does not matter in 3.1 if we demand ``$\beta^\delta_n$ is a successor
ordinal".
\endproclaim
\bigskip

\demo{Proof}  0) Check. \newline
1)  By 2). \newline
2), 3)  Check. \hfill$\square_{3.2}$
\enddemo
\bigskip

\subhead {Questions} \endsubhead  1) Are there such examples for
$\lambda$ singular? \newline
2)  Suppose in the definition we allow for each $\alpha$ a filter on
$Dom(f_\alpha)$ generated by $\aleph_0$ sets; do we get an equivalent
principle?
\bigskip

\proclaim {3.3 Claim}  Assume $INCwh^2(\kappa)$, $\lambda > \kappa$,
$\lambda = cf(\lambda)$,
$S \subseteq \{ \delta < \lambda:cf(\delta) = \aleph_0\}$ is stationary and
 $\diamondsuit_S$ holds. \newline
\underbar{Then} (1) there is a $\langle \langle m^\delta_n,\beta^\delta_n:n <
\omega \rangle: \delta \in S \rangle$ as in 2.4(2), but only (a) and
$(b)_\kappa$ hold.
\newline
2)  For some regular $\lambda' \in [\kappa,\lambda]$, we have $INCWH^2
(\lambda')$. \newline
3)  We can replace $INCwh^2(\lambda')$, $INCWH^2(\lambda)$ by
$INCwh^1(\lambda)$, $INCWH^1(\lambda')$ respectively.
\endproclaim
\bigskip

\demo{Proof}  Now (2) follows from (1) as in 2.7 and we leave (3) to the
reader.  The proof of 3.3(1) is like the
proof of 1.3 with one twist.  Let $h:\lambda \rightarrow \kappa$ be such
that for every $\zeta < \kappa$, $h^{-1}(\{ \zeta \})$ has cardinality
$\lambda$.  Let
$\left\langle \langle m^\delta_n,{}^* \beta^\delta_n:n <
\omega \rangle:\delta \in S^* \right\rangle$ \newline
witness $INCwh^2(\kappa)$.

Let $\langle g_\delta:\delta \in S \rangle$ witnesses $\diamondsuit_S$ i.e.
$g_\delta \in {}^\delta \omega$ and for every $g \in {}^\lambda \omega$ for
stationarily many $\delta \in S$, $g_\delta = g\restriction \delta$. \newline
For each $\delta \in S$ we define a function $g^*_\delta \in {}^\kappa
\omega$:

$$
\align
g^*_\delta(\zeta) = \text{Min} \biggl\{ m:&\text{for arbitrarily large }
\alpha < \delta \text{ we have}:m = g_\delta(\alpha) \text{ and } \\
  &h(\alpha) = \zeta \biggr\}, \text{ if defined}.
\endalign
$$

If for some $\zeta < \kappa$, $g^*_\delta(\zeta)$ is not defined
(i.e. there is no such  $m$) - we do nothing.

If $g^*_\delta \in {}^\kappa \omega$ is defined we know that for some
$\zeta(\delta) \in S^*$, $(\exists^\infty n)(g^*_\delta({}^*\beta
^{\zeta(\delta)}_n) \le m^\delta_n)$.  
\newline
(Such a $\zeta(\delta)$ exists by the choice of
$\left\langle \langle m^\zeta_n,{}^*\beta^\zeta_n:
n < \omega \rangle:\zeta \in S^* \right\rangle)$.  We fix such a $\delta$.
\newline
For each $n < \omega$ choose $\xi(\delta,n) < \kappa$ such that:
\medskip
\roster
\item "{$(*)_1$}"  for arbitrarily large $\gamma < \delta$ we have \newline
$(*)_{\gamma,n} \, \gamma \notin S,g_\delta(\gamma) =
 g^*_\delta({}^*\beta^{\zeta(\delta)}_n) \wedge h_0(\gamma) = {}^*\beta
^{\zeta(\delta)}_n \wedge h_1(\gamma) = \xi(\delta,n)$
\item "{$(*)_2$}"  $\zeta = \zeta(\delta) < \xi(\delta,n)$.
\endroster
\medskip

\noindent
Choose $\gamma^\delta_n$ such that:
\medskip
\roster
\item "{(a)}"  $\gamma^\delta_n < \delta,
h(\gamma^\delta_n) = {}^*\beta^{\zeta(\delta)}_n,g_\delta
(\gamma^\delta_n) = g^*_\delta({}^*\beta^{\zeta(\delta)}_n)$
\item "{(b)}"  $\delta = \dsize \bigcup_{n < \omega} \gamma^\delta_n$ and
$\gamma^\delta_n < \gamma^\delta_{n+1}$.
\endroster
\medskip

\noindent
We claim $\left\langle \langle m^{\zeta(\delta)}_n,\gamma^\delta_n:n < \omega
\rangle:\delta \in S \right\rangle$ witness the conclusion.  Looking at
Definition 2.4, the preliminary properties hold. \newline
We have to prove clause $(a)$ of 2.7.
\enddemo
\bigskip

\demo{Proof of (a)}  Let $g \in {}^\lambda \omega$.  For each $\zeta <
\kappa$, $\{ \alpha < \lambda:h(\alpha) = \zeta \}$
has cardinality $\lambda$, so

$$
g^*(\zeta) = \text{Min}\{ m:(\exists^\lambda \alpha)[h(\alpha) = \zeta
\wedge g(\alpha) = m] \}
$$

\noindent
is well defined.  Let

$$
A =: \{ (\zeta,m):(\exists^\lambda \alpha < \lambda)[g(\alpha) = m
 \wedge h(\alpha) = \zeta] \text{ and } \zeta < \kappa,m < \omega \}.
$$

\noindent
Then

$$
\align
E =: \biggl\{ \delta < \lambda:&\text{ for every } (\zeta,m) \in A,
\text{ for } \lambda \text{ many } \\
  &\alpha < \lambda, g(\alpha) = m,h(\alpha) = \zeta,\text{ and for every }\\
  &(\zeta,m) \in (\kappa \times \omega) \backslash A, \text{ we have } \\
  &\delta > \text{ sup}\{ \alpha < \lambda:g(\alpha) = g^*(\zeta) \wedge h
(\alpha) = \zeta \} \biggr\}
\endalign
$$

\noindent
is a club of $\lambda$. \newline
For stationarily many $\delta \in S$, $g_\delta \subseteq g$ so there is 
such a $\delta \in E \cap S$.

Now check: $g^*_\delta = g^*(g^*_\delta$ was defined earlier).
The rest is also easy to check.
\enddemo
\bigskip

\demo{Proof of $(b)_\lambda$ i.e. $(< \kappa)$-freeness}
 Let $u \subseteq \lambda$, $|u| <
\kappa$, hence $v = \{ h(\alpha):\alpha \in u \}$ is a subset of $\kappa$
of cardinality $< \kappa$, so by the choice of $\langle m^\delta_n,
{}^* \beta^\delta_n:n < \omega,\delta \in S \rangle$ there is a 
$f^*:v \rightarrow \omega$ as required. \newline
Choose $f:u \rightarrow \omega$ by $f(\alpha) = f^*(h(\alpha))$, now $f$ is
as required. \hfill$\square_{3.3}$
\enddemo
\bigskip

\demo{3.4 Discussion}  1) Probably $INCWH(\lambda)$ should mean just there is
a first countable ($<\lambda$)-CWH not $\lambda$-CWH, as this is actually
the two notions which speak on $m_\alpha,\beta^\delta_n(n < \omega)$ or
$m^\alpha_n,\beta^\delta_n(n < \omega)$ and they should be named
$INCWH^\ell(\lambda)$,$\ell = 1,2$, respectively. \newline
So, $(\exists \lambda \ge \mu)INCWH {}^\ell (\lambda)$ is equivalent to
$(\exists \lambda \ge \mu)INCwh {}^\ell (\lambda)$ (for $\ell = 1,2$).
\enddemo
\bigskip

\definition{3.5 Definition}  1)  $INCWH^3(\lambda)$ means: there are
$S \subseteq \lambda$ and $f:\lambda \times \lambda \rightarrow \omega$
such that if we define the spaces as before, i.e.

$$
\text{the points of } X \text{ are } 
y_\alpha, x_{\alpha,\beta}, (\alpha < \beta < \lambda)
$$

$$
\text{each } x_{\alpha,\beta} \text{ is isolated}
$$

$$ 
\align
u_{\alpha,n} = \{ y_\alpha \} &\cup \{x_{\alpha,\beta}:f
(\alpha,\beta) \le n,\alpha < \beta,\alpha \notin S,\beta \in S \} \\
   &\cup \{ x_{\beta,\alpha}:f(\beta,\alpha) \le n,
\beta < \alpha,\beta \notin S,\alpha \in S \}
\endalign
$$

\noindent
for $n \in \omega$ is a neighborhood base at $Y_\alpha$ such that:
\medskip
\roster
\item "{(a)}"  $\alpha < \beta < \lambda$, $u_{\alpha,n} \cap u_{\beta,m} \ne
\emptyset \Rightarrow \beta \in S \wedge \alpha \notin S$
\item "{(b)}"  for every $\alpha < \beta \in S$ and for some $n$ we have:
$\alpha < \gamma < \beta \Rightarrow u_{\beta,n} \cap u_{\gamma,0} =
\emptyset$,
\endroster
\noindent
then
\roster
\item "{(c)}"  the space $X$ is not CWH but is ($< \lambda$)-CWH.
\endroster
\medskip
\noindent
2)  $INCWH^4(\lambda)$ means: there is a symmetric two-place function $f$
from $\lambda$ to \newline
$\{ v:v \subseteq \omega \times \omega$ is finite, and $(n,m) \in v,n' \le
n,m' \le m \Rightarrow (n',m') \in v \}$ which is not free (i.e. for any
$g:\lambda \rightarrow \omega$ for some $\alpha < \beta$,
$(g(\alpha),g(\beta)) \in f(\alpha,\beta)$), but is $\lambda$-free (i.e.
for every $A \subseteq \lambda$, $|A| < \lambda$, there is a $g:A \rightarrow
\omega$ with no such $\alpha < \beta$ which are from $A$.
\enddefinition
\bigskip

\noindent
The point is that
\proclaim{3.6 Claim}  1)  $INCWH^1(\lambda) \Rightarrow INCW^2 (\lambda)
\Rightarrow INCWH^3(\lambda) \Rightarrow INCWH^4(\lambda)$. \newline
2)  \underbar{If} $\lambda = cf(\lambda) > \aleph_0$ and
$S \subseteq \{ \delta < \lambda:cf(\delta) = \aleph_0 \}$ is stationary
not reflecting, \underbar{then} $INCWH^3(\lambda)$.
\endproclaim
\bigskip

\proclaim{3.7 Lemma}  In 2.5 we can weaken $INCWH^1_\lambda$ to
$INCWH^3(\lambda)$.
\endproclaim
\bigskip

\demo{Comment}  The $INCWH^\ell(\lambda)$ are not so artificial: we can
translate \newline
$INCWH(\lambda)$ to a similar statement.
\enddemo
\bigskip

\proclaim{3.8 Claim}  $SINCWH(\lambda) \Rightarrow INCWH^4(\lambda)$.
\endproclaim
\bigskip

\demo{Proof}  Let the space $X$ exemplify $SINCWH(\lambda)$.
Let $\{ y_\alpha:\alpha < \lambda \} \subseteq X$ exemplifies 
``$X$ not $\lambda-CWH$" i.e. it is discrete not separated and
$\alpha \ne \beta \Rightarrow y_\alpha \ne y_\beta$.

Let $u_{\alpha,n} \supseteq u_{\alpha,n+1}, \{ u_{\alpha,n}:n < \omega
\}$ be a neighborhood basis of $y_\alpha$.  Now for each $\alpha,n,\beta,m$
choose
if possible $x_{\alpha,n,\beta,m} \in u_{\alpha,n} \cap u_{\beta,m}$.
Let  $f(\alpha,\beta) = \{ (n,m):x_{\alpha,n,\beta,m} \text{ is defined}\}$.
This $f$ exmplifies $INCWH^4(r)$. \hfill$\square_{3.8}$
\enddemo
\bigskip

\remark{3.8A Remark}  The $\Leftarrow$ holds as well.
\endremark 
\newpage

\head {\S4 General Set Theoretic Spectrum of Freeness} \endhead
\bigskip

\definition{4.0 Definition}  For $\lambda > cf(\lambda) = \theta$ let
$(*)_\lambda$ means: there is a
$\{ \eta_\alpha:\alpha < \lambda^+ \} \subseteq
{}^\theta \lambda$ which is $\lambda$-free (see (c) in 4.1(1) below).
\enddefinition
\bigskip

\definition{4.1 Definition}  1) For $\theta$ a regular cardinal and
$\sigma \ge 1$ (if $\sigma = 1$ we omit it) let:
$$
\align
SP_{\theta,\sigma} = \biggl\{ \lambda:&\text{there is a family } H
\text{ such that}: \\
  &(a) \, \text{ every } h \in H  \text{ is a partial function from
ordinals to } \theta \\
  &(b) \,\, h \in H \Rightarrow |Dom(h)| = \theta \\
  &(c) \, \text{every } H' \subseteq H \text{ of cardinality } < \lambda
 \text{ is } \sigma \text{-free which means that } \\
&\quad \text{ it can be represented as a union } 
\dsize \bigcup_{i<i(*)} H'_i,i(*) < 1 + \sigma, \\
&\quad \text{ and each } H'_i \text{ is free.  For } H'_i \text{ to be
free means that there} \\
&\quad \, \text{is a } g, \text{ a function from ordinals to }
 \theta \text{ such that } \\
&\quad \, (\forall h)(\exists \xi < \theta)[h \in H'_i \rightarrow \\
&\quad \, (\forall \alpha \in \text{ Dom}(h)[h(\alpha) \le g(\alpha) \vee
 h(\alpha) \le \xi] \\
&(d) \,\, H \text{ is not } \sigma \text{-free, } |H| = \lambda \biggr\}
\endalign
$$
\medskip

\noindent
2)
$$
\align
SPd_{\theta,\sigma} = \{ \lambda:&\text{ there is an } H 
\text{ satisfying (a)-(d) above and} \\
  &\text{(e) each } h \in H \text{ is one to one} \}.
\endalign
$$
\newpage

3)
$$
\align
SPw_{\theta,\sigma} = \biggl\{ \lambda:&\text{there is a family }
H \text{ such that}: \\
&(a) \, \text{ if } (h,\bar u) \in H  \text{ then } h
\text{ is a function from ordinals to } \theta \\
&(b) \, \text{ if } (h,\bar u) \in H, \text{ then }
\bar u = \langle u_\varepsilon:\varepsilon < \theta \rangle \text{ is a
decreasing sequence} \\
&\quad \, \text{ of subsets of } Dom(h) \\
&(c)  \text{ every pain } (H^1,Z^1),
\text{ with } Z^1 \le \text{ ordinals}, \\
&\quad \, |Z^1| < r \text{ and } H' \subseteq H  \text{ of cardinality }
< \lambda \text{ is } \sigma\text{-free, which means} \\
&\quad \, \text{ it can be represented as } \dsize \bigcup_{i<i(*)}
(H^1_i,Z^1_i),i(*) < 1 + \sigma \\
&\quad \, \text{ and each } (H'_i,Z'_i) \text{ is free.  This means that
there are functions} \\
&\quad \, g,f \text{ with } g:H'_i \rightarrow \theta \text{ and } f
\text{ from ordinals to } \theta \\
&\quad \, \text{ such that for every } (h,\bar u) \in H', \text{ for every }\\
&\quad \, \alpha \in Z'_i \cap \, Dom(h) \text{ we have } h(\alpha) \le
\, max\{f(\alpha),g(h)\}. \\
&(d) \,  H \text{ is not } \sigma\text{-free},  |H| = \lambda \biggr \}
\endalign
$$
\enddefinition
\bigskip

\demo{4.2 Observation}  0) In Definition 4.1, if each $h \in H$ converges
to $\theta$, in clause (c) of 4.1(1) we can just demand
$(\forall h)[h \in H' \rightarrow \theta > | \{ \alpha:
h(\alpha) > g(\alpha) \}|]$. \newline
1)  In Definition 4.1(1) without loss of generality $\dsize \bigcup
_{h \in H} \text{Dom}(h) \subseteq \lambda$ and in 4.1(3) without loss of
generality $\dsize \bigcup \Sb (h(\bar u) \in H) \\ (h,\bar u) \endSb
\text{Dom}(h) \subseteq \lambda$.  Also, without loss of generality
$\text{Dom}(g) = \lambda$. \newline
2)  Note  \newline
$\theta^+ \cap SP_\theta = \emptyset$
[why? if $H = \{ h_\zeta:\zeta < \zeta^* \le \theta \}$,
$\dsize \bigcup_{\zeta} \text{ Dom}(h_\zeta) = \{ \alpha_i:i <
\theta\}$, let $g(\alpha_i) = \text{ sup}\{h_\zeta(\alpha_i):\zeta < i, \,
\alpha_i \in \text{ Dom}(h_\zeta) \}$].  This also follows from 4.1(B1) and
4.2(2).  \newline
3)  $SP_\theta \cap [\theta^+,2^\theta] \ne \emptyset$ [this follows from
4.2(4) below]. \newline
\newpage

4) We let
$$
\align
{\frak b}[\theta] = \text{ Min} \{|F|: \,&F \subseteq
 {}^\theta \theta,\text{ and for no } g \in
 {}^\theta \theta \text{ do we have }\\
  &(\forall f \in F)(\exists \zeta < \theta)(f\restriction[\zeta,\theta)
< g \restriction [\zeta,\theta)) \}
\endalign
$$
\noindent
if $\sigma \le \theta^+$ then clearly ${\frak b}[\theta] \in SP_{\theta,
\sigma}$. \newline
5)  In Definition 4.1(3) without loss of generality for $(h,\bar u) \in H$,
$\dsize \bigcap_{\zeta < \theta} u_\zeta = \emptyset$.  Also without loss
of generality, for $(k,\bar u) \in H, u_\zeta = \{ \alpha \in \text{ Dom}(h):
h(\alpha \ge \zeta \}$ (we say: $\,\bar u$ is standard for $h$). \newline
6)  Suppose that $H$ is as in 4.1(3). In c), if we set $Z' = \theta$ and
assume that $\bar u$ is standard, we obtain: \newline
For every $H' \subseteq H$ with $|H'| < \lambda$, there are sets $H'_i$ for
$i < i(*) < 1 + \sigma$ such that $H' = \dsize \bigcup_{i < i(*)} H'_i$ and
for each $i < i(*)$, there is a function $g_i:H'_i \rightarrow \theta$ with
the following property. \newline
For every $(k,\bar u) \in H'_i$ \newline
$\exists \xi < \theta,\exists \xi < \theta,\forall \alpha \in u_\xi[k(\alpha)
\le \text{ max}\{\xi,g_i(\alpha) \}]$. \newline
7)  Note also that we can without loss of generality assume that
$Z' \subseteq \dsize \bigcup_{k \in H'} \text{ Dom}(k)$, for
4.1.3)c). \newline
8)  We restrict our attention to the case $\sigma \le \theta^+$.  Actually,
the main interest is in $\sigma = 1$.  For $\sigma$ large enough the
definition of $\sigma$-free sets as it stands would imply that all relevant
$H$ are $\sigma$-free, if $|H| = \lambda$.
\enddemo
\bigskip

\subhead {Notation} \endsubhead  For $a \subseteq \theta \times \theta:a$ is
pic if $(\zeta_1,\xi_1) \ne (\zeta_2,\xi_2) \in a \Rightarrow \neg
(\zeta_1,\xi_1) \le (\zeta_2,\xi_2)$ coordinatewise. \newline
Pic$(\theta \times \theta) = \{a:a \subseteq \theta \times \theta$ and $a$
is pic (hence finite)$\}$ \newline
$Cl(a) = \{ (\zeta,\xi) \in (\theta \times \theta):(\exists x \in a)
(x \le (\zeta,\xi)$ coordinatewise$\}$, for $a \subseteq \theta \times
\theta$.
\newpage

\definition{4.1A Definition}  1) For $\theta$ a regular cardinal and
$\sigma \ge 1$ (if $\sigma = 1$ we omit it) let:
$$
\align
SQ_{\theta,\sigma} = \biggl\{ \lambda:&\text{there is a family } H
\text{ such that}: \\
  &(a) \, \text{ every } h \in H  \text{ is a partial function from the
ordinals } \\
  &(b) \,\, h \in H \Rightarrow |Dom(h)| = \theta,Rang(h) \subseteq Pic
(\theta \times \theta) \\
  &(c) \, \text{every } H' \subseteq H \text{ of cardinality } < \lambda
 \text{ is } \sigma \text{-free which means that } \\
&\quad \text{ it can be represented as a union } 
\dsize \bigcup_{i<i(*)} H'_i,i(*) < 1 + \sigma, \\
&\quad \text{ and each } H'_i \text{ is free.  For } H'_i \text{ to be
free means that there} \\
&\quad \, \text{is a } g, \text{ a function from ordinals to }
 \theta \text{ such that } \\
&\quad \, (\forall h)(\exists \xi < \theta)[h \in H'_i \rightarrow
(\forall \alpha \in \text{ Dom}(h)[(g(\alpha),\Xi) \in Cl(k(\alpha))] \\
&(d) \,\, H \text{ is not } \sigma \text{-free, } |H| = \lambda \biggr\}.
\endalign
$$
\medskip

\noindent
2)
$$
\align
SQd_{\theta,\sigma} = \biggl\{ \lambda:&\text{ there is an } H 
\text{ satisfying (a)-(d) above and} \\
  &\text{(e) each } h \in H \text{ is simple, which means: there is an} \\
  &\quad \, \text{ enumeration } Dom(h) = \{ \alpha_\zeta:\zeta < \theta
\} \text{ with no repetitions}, \\
  &\quad \, \text{ such that } h(\alpha_g) = \{(\zeta_1,\zeta_2):(\zeta_1,
\zeta_2) \nleq (\beta_\zeta,\gamma_\zeta)\} \\
  &\quad \, \text{ for some } \langle \gamma_\zeta:\zeta < \theta \rangle
\text{ which are strictly increasing and} \\
  &\quad \, \dsize \bigcup_{\xi < \zeta} \beta_\zeta < \gamma_\zeta \biggr\}.
\endalign
$$
\newpage

3)
$$
\align
SQw_{\theta,\sigma} = \biggl\{ \lambda:&\text{there is a family }
H \text{ such that}: \\
&(a) \, \text{ if } (h,\bar u) \in H  \text{ then } h
\text{ is a function from ordinals to } Pie(\theta \times \theta) \\
&(b) \, \text{ if } (h,\bar u) \in H, \text{ then }
\bar u = \langle u_\varepsilon:\varepsilon < \theta \rangle \text{ is a
decreasing sequence} \\
&\quad \, \text{ of subsets of } Dom(h) \\
&(c)  \text{ every pain } (H^1,Z^1),
\text{ with } Z^1 \le \text{ ordinals}, \\
&\quad \, |Z^1| < r \text{ and } H' \subseteq H  \text{ of cardinality }
< \lambda \text{ is } \sigma\text{-free, which means} \\
&\quad \, \text{ it can be represented as } \dsize \bigcup_{i<i(*)}
(H'_i,Z'_i),i(*) < 1 + \sigma \\
&\quad \, \text{ and each } (H'_i,Z'_i) \text{ is free.  This means that
there are functions} \\
&\quad \, g,f \text{ with } g:H'_i \rightarrow \theta \text{ and } f
\text{ from ordinals to } \theta \\
&\quad \, \text{ such that for every } (h,\bar u) \in H'_i,
\text{ for every }\\
&\quad \, \bar z \in Z'_i \cap \, Dom(h) \text{ we have } \\
&\quad \, (g(h),f(z)) \in c\ell(k(z)) \\
&(d) \,  H \text{ is not } \sigma\text{-free},  |H| = \lambda \\
&(e) \, (k,\bar u) \in H \Rightarrow \dsize \bigcap_{\varepsilon < \theta}
u_\varepsilon = \emptyset \biggr \}.
\endalign
$$

\noindent
\underbar{Note}:  1) In 4.1A3)c), we can assume that $Z' \subseteq
\dsize \bigcup_{h \in H'} Dom(h)$. \newline
2)  As in 4.1, we consider only the case $\sigma \le \theta^+$. \newline
3)  $SPx_{\theta,\sigma}$ can be understood as a particular case of
$SQx_{\theta,\sigma}$, where $\text{Rang}(h)$ is restricted to
$\{ (\zeta,\zeta):\zeta < \theta \}$.  Here, $x \in \{ w,d \}$ or $x$ is
omitted.

\enddefinition
\bigskip

\demo{4.1B Fact}  1) $\lambda \in SP_{\theta,\sigma}$ implies that
$\lambda \in SQ_{\theta,\sigma}$ \newline
$\lambda \in SPd_{\theta,\sigma}$ implies that $\lambda \in
SQd_{\theta,\sigma}$, and \newline
$\lambda \in SPw_{\theta,\sigma}$ implies that
$\lambda \in SQw_{\theta,\sigma}$. \newline
2)  $\lambda \in SQd_{\theta,\sigma}$ implies that
$\lambda \in SP_{\theta,\sigma}$.
\enddemo
\bigskip

\demo{Proof}  $H$ exmplifies that $\lambda \in SP_{\theta,\sigma}$, let
$H^\otimes = \{ h^x:h \in H \}$, where for $h \in H$, $h^\otimes$ is a
function with domain $Dom(h)$ and
$$
h^\otimes(\alpha) = \{ (h(\alpha),h(\alpha))\}.
$$


\noindent
Similarly for $SPd_{\theta,\sigma}$. \newline
If $H$ exmplifies that $\lambda \in SQw_{\theta,\sigma}$, let $H^\otimes =
\{ (h^\otimes,\bar u):(h,\bar u) \in H \}$. \newline
2)  Let $H = \{ h_j:j < \lambda \}$ exemplifies that $\lambda \in
SQd_{\theta,\sigma}$,  Let us enumerate $Dom(h_j) \le \{ \alpha^j_\zeta:
\zeta < \theta \}$ for $j < \lambda$, as in clause (e) of 4.1A(2). \newline
Then we know that

$$
\align
h_j(\alpha^j_\zeta) = \biggl\{ (\varepsilon_1,\varepsilon_2):&\varepsilon_1
< \theta \text{ and } \varepsilon_2 < \theta \text{ and } \\
  &(\varepsilon_1,\varepsilon_2) \nleq (\beta^j_\zeta,\gamma^j_\zeta)
\biggr\}
\endalign
$$

\noindent
for some $(\gamma^j_\zeta:\zeta < \theta)$ which is strictly increasing
and $\gamma^j_\zeta > \dsize \bigcup_{\xi < \zeta} \beta^j_\zeta$. \newline
Let $h^\oplus_j$ be the function with domain $Dom(h'_j) = \{ \alpha^j
_\zeta:\zeta < \theta \}$ and defined by $h^\oplus_j(\alpha_\zeta) = 
\beta^j_\zeta$.  Then $H^\oplus = \{ h^\oplus_j:j < \lambda \}$ exemplifies
that $\lambda \in SP_{\theta,\sigma}$.
\enddemo
\bigskip

\subhead {Notation} \endsubhead  For a function $h$ from a subset of
ordinals to $Pie(\theta \times \theta)$, we say that $h$ converges to
$\theta$, if
$$
\align
(\forall \beta < \theta)&(\exists \alpha)(\forall \gamma \in \text{ Dom}(h)
\backslash \alpha) \\
  &[(\varepsilon_1,\varepsilon_2) \in h(\gamma) \Rightarrow \varepsilon_1
> \beta \text{ and } \varepsilon_2 > \beta].
\endalign
$$

\demo{4.2A Observation}  0) In Definition 4.1.A, if each $h \in H$ converges
to $\theta$, in clause (c) of 4.1A(2) we can just demand \newline
$(\forall h)[h \in H' \rightarrow \theta > | \{ \alpha:
\exists(\varepsilon_1,\varepsilon_2) \in k(\alpha)[\varepsilon_1 > g(\alpha)
\vee \varepsilon_2 > g(\alpha)]|$. \newline
1)  In Definition 4.1(1) without loss of generality $\dsize \bigcup
_{h \in H} \text{Dom}(h) \subseteq \lambda$ and in 4.1A(3) without loss of
generality $\dsize \bigcup \Sb (h(\bar u) \in H) \\ (h,\bar u) \endSb
\text{Dom}(h) \subseteq \lambda$.  Also, without loss of generality,
$\text{Dom}(g) = \lambda$. \newline
2)  Note  \newline
$\theta^+ \cap SQ_\theta = \emptyset$
[why? if $H = \{ h_\zeta:\zeta < \zeta^* \le \theta \}$,
$\dsize \bigcup_\zeta \text{ Dom}(h_\zeta) = \{ \alpha_i:i <
\theta\}$, let $g(\alpha_i) = \text{ sup} \{\text{max }h_\zeta(\alpha_i):
\zeta < i, \, \alpha_i \in \text{ Dom}(h_\zeta) \}$]. \newline
3)  $SQ_\theta \cap [\theta^+,2^\theta] \ne \emptyset$ [this follows from
4.2.3) and 4.1.B1)].  Actually, $b[\theta] \in SQ_\theta$.
4)  In 4.1A3)c), if we set $Z' = \theta$, we obtain the following property.
\newline
For every $H' \subseteq H$ of cardinality $< \lambda$, there are sets
$H'_i$ for $i < i(*) < 1 + \sigma$, such that there are functions
$(q_i:i < i(*)),q_i:H'_i \rightarrow \theta$ satisfying: if $(h,\bar u) \in
H'_i$, then \newline
$(\exists \zeta < \theta)(\exists \xi < \theta)(\forall \alpha \in u_\zeta)
[(q_i(\alpha),\xi) \in c \ell(h(\alpha)]$.
\enddemo
\bigskip

\proclaim{4.3 Claim}  1) If there is an $H$ as in (a), (b) of 4.1(1) which
is $(< \mu)-\sigma$-free not $\lambda-\sigma$-free \underbar{then} there
is a $\lambda' \in [\mu,\lambda] \cap SP_{\theta,\sigma}$.  Similarly for
4.1(2), 4.1(3). \newline
2)  If $pp_{\varGamma(\theta)}(\lambda) > \lambda^+$, $\lambda > cf(\lambda)
= \theta$ (or just $(*)_\lambda$ of 4.0) and \newline
$\lambda \ge \sigma$
\underbar{then} $SP_{\theta,\sigma} \cap
 [\lambda^+,\lambda^\theta] \ne \emptyset$.
\endproclaim
\bigskip

\demo{Proof}  1) Straightforward. \newline
2)  Let for $\{ \eta_\alpha:\alpha < \lambda^+ \}
 \subseteq {}^\theta \lambda$ be $\lambda$-free, without loss of generality
\newline
$\langle \{ \eta_\alpha(\zeta):\alpha < \lambda^+ \}:\zeta < \theta \rangle$
are pairwise disjoint and let

$$
\align
H = \biggl\{ h:\,&\text{for some } \alpha < \lambda^+ \text{ and } a \subseteq
\lambda^+, \text{otp}(a) = \theta,\text{Dom}(h) = a,  \\
 &h \text{ is strictly increasing and for } \beta \in a \\  
&h(\beta) = \text{ sup}\{ \varepsilon:\eta_\alpha(\varepsilon) = \eta_\beta
(\varepsilon) \} \biggr\}.
\endalign
$$
\medskip

\noindent
Now $H$ is not free: if $g:\lambda^+ \rightarrow \theta$, then for some
$\varepsilon < \theta$, $A = \{ \alpha < \lambda^+:g(\alpha) = \varepsilon
\}$ is of cardinality $\lambda^+$.  Choose by induction on
$\zeta < \lambda^+$ an ordinal $\alpha^*_\zeta < \lambda^+$
increasing with $\zeta$ such that \newline
$\cup \{ \text{Rang}(\eta_\alpha):\alpha \in A \cap \alpha^*_{\zeta+1}
\backslash \alpha^*_\zeta \} = \cup \{ \text{Rang}(\eta_\alpha):\alpha \in A
\backslash \alpha^*_\zeta \}$. \newline
Next choose $\alpha \in A \backslash \alpha^*_\theta$ and
$\beta_\zeta \in A \cap [\alpha^*_\zeta,\alpha^*_{\zeta+1})$ such that
$\eta_{\beta_\zeta}
(\zeta) = \eta_\alpha(\zeta)$ and let $a = \{ \beta_\zeta:\zeta < \theta \}$,
$h \in {}^a \theta$, $h(\beta_\zeta) = 
\text{ sup}\{ \varepsilon:\eta_\alpha(\varepsilon) = \eta_{\beta_\zeta}
(\varepsilon)\} \ge \zeta$, so $h \in H$.  As
$\beta_\zeta \in A$, $g(\beta_\zeta) = \varepsilon=$constant, so if
$\xi < \theta$, $\{ \beta \in \text{ Dom}(h):h(\beta) \ge g(\beta),
\xi \}$ include $\{ \beta_\zeta:\xi,\varepsilon < \zeta < \theta \}$, which
is a contradiction.

On the other hand, $H$ is $\lambda^+$-free.  For suppose $H' \subseteq H$,
$|H'| \le \lambda$.  For $h \in H'$ choose $\alpha_h,a_h$ witnessing
$h \in H$.  Then $b = \cup \{ \{ \alpha_h \} \cup a_h: h \in H' \}$
is a subset of $\lambda^+$ of cardinality $\le \lambda$, hence we can find
$\langle \varepsilon_\alpha:\alpha \in b \rangle$ such that
$\langle \text{Rang}(\eta_\alpha \restriction [ \varepsilon_\alpha,
\theta)):\alpha \in b \rangle$ is a sequence of pairwise disjoint subsets
of $\lambda$.  Let us define a $g:\lambda^+ \rightarrow \theta$ such that
$\alpha \in b \Rightarrow g(\alpha) = \varepsilon_\alpha$.  Now if
$h \in H'$, let $a_h = \{ \beta_\zeta:\zeta < \theta \}$ (increasing with
$\zeta$), so

$$
\align
h(\beta_\zeta) = &\text{ sup}\{ \varepsilon:\eta_{\alpha_h}(\varepsilon)
= \eta_{\beta_\zeta}(\varepsilon) \} \\
  &\text{ so } h(\beta_\zeta) \le \text{ max}\{ \varepsilon_{\alpha_h},
\varepsilon_{\beta_\zeta} \} = \text{ max}\{ g(\alpha_h),g(\beta_\zeta) \}
\endalign
$$

\noindent
So choose $\xi = g(\alpha_h)$ and we get the desired conclusion. \newline
To finish we use part (1).\hfill$\square_{4.3}$
\enddemo
\bigskip

\proclaim{4.3A Claim}  1) If there is an $H$ as in (a), (b) of 4.1A(1) which
is $(< \mu)-\sigma$-free not $\lambda-\sigma$-free \underbar{then} there
is a $\lambda' \in [\mu,\lambda] \cap SQ_{\theta,\sigma}$.  Similarly for
4.1A(2), 4.1A(3). \newline
2)  If $pp_{\varGamma(\theta)}(\lambda) > \lambda^+$, $\lambda > cf(\lambda)
= \theta$ (or just $(*)_\lambda$ of 4.0) and \newline
$\lambda \ge \sigma$
\underbar{then} $SQ_{\theta,\sigma} \cap
 [\lambda^+,\lambda^\theta] \ne \emptyset$.
\endproclaim
\bigskip

\demo{Proof}  1) Straightforward. \newline
2) This follows from 4.3.2) and 4.1.B1). \hfill$\square_{4.3A}$
\enddemo
\bigskip

\proclaim{4.4 Claim}  1)  The following implications hold for any $\lambda$:
$$
(a) \Rightarrow (b) \Leftrightarrow (b)^+ \Leftrightarrow (c)
\Leftrightarrow (c)^+  \Rightarrow (d),
$$

\noindent
where
\roster
\item "{$(a)$}"  $\lambda \in SQ_{\lambda_0}$
\item "{$(b)$}"  There is a $(< \lambda)$-CWH not $\lambda$-CWH
first countable space
\item "{$(b)^+$}" There is a space like in (b), which is in addition
$(< \lambda)$-metrizable
\item "{$(c)$}"  There is a $(< \lambda)-^*$ CWN first countable space
with $\lambda$ points.
\item "{$(c)^+$}"  Thre is a space like in (c), which is in addition
$(< \lambda)$-metrizable
\item "{$(d)$}"  $\lambda \in SQw_{\aleph_0}$.
\endroster
\noindent
2) $\lambda \in SQ_{\theta,\sigma} \Rightarrow \lambda \in 
SQw_{\theta,\sigma} \Rightarrow [\lambda,\lambda^\theta] \cap
SQd_{\theta,\sigma} \ne \emptyset$ for $\sigma \le \theta^+$. \newline
3) $\lambda \in SP_{\theta,\sigma} \Rightarrow \lambda \in 
SPw_{\theta,\sigma} \Rightarrow [\lambda,\lambda^\theta] \cap
SPd_{\theta,\sigma} \ne \emptyset$ for $\sigma \le \theta^+$. \newline
4)  Similarly for ${\Cal T}_\theta$.
\endproclaim
\bigskip

\demo{Proof}  1),4)(a) implies (b), (b)$^+$, (c), (c)$^+$.
First implication - assume $H$ exemplifies that $\lambda \in SQ_\theta$,
we can use the space \newline
$X = \{ y_i:i <
\lambda \} \cup \{ z_h:h \in H \} \cup \{ x_{h,i}:h \in H,i \in \text{ Dom}
(h) \}$, and for \newline
$\zeta < \theta$ let $u_\zeta[z_h] = \{z_h \} \cup
\{x_{h,i}:i \in \text{ Dom}(h),(\zeta,\zeta) \notin c \ell(h(i))\}$, \newline 
$u_\zeta[y_i] =
\{y_i \} \cup \{ x_{h,i}:h \in H,i \in \text{ Dom}(h),(\zeta,\zeta) \notin
c \ell(y_i(i)) \}$ and $x_{h,i}$ is isolated.

Suppose $H' \subseteq H$, $|H'| < \lambda$ and let \newline
$X[H'] = \{ y_i:i < \lambda \}
\cup \{ z_h:h \in H' \} \cup \{ x_{h,i}:h \in H, i \in \text{ Dom}(h) \}$.
\newline
Let $g:\lambda \rightarrow \theta$ be such that for every $h \in H'$, for some
$\zeta[h] < \theta$ we have

$$
i \in \text{ Dom}(h) \Rightarrow g(i),\zeta[h]) \in c \ell(h(i)).
$$

\noindent
Let us choose for $t \in X[H']$ a neighborhood $v_t$:

$$
\alignat2
\text{\underbar{if}} \qquad t &= x_{h,i} &&\qquad
\text{\underbar{then}} \qquad v_t = \{x_{h,i} \} \\
\text{\underbar{if}} \qquad t &= y_i &&\qquad \text{\underbar{then}} \qquad
v_t = u_{g(i)}[y_i] \\
\text{\underbar{if}} \qquad t &= z_h &&\qquad \text{\underbar{then}} \qquad
v_t = u_{\zeta[h]}[z_h].
\endalignat
$$

\noindent
Now

$$
\langle v_{y_i}:i < \lambda \rangle \char 94 \langle v_{z_h}:h \in H'
\rangle \char 94 \langle v_{x_{h,i}}:i < \lambda,h \in H \text{ and }
x_{h,i} \notin \dsize \bigcup_{j<\lambda} v_{y_j} \cup \dsize
\bigcup_{h \in H'} v_h \rangle
$$

\noindent
is a partition of $X[H']$ to pairwise
disjoint open sets.  In each basic open set there is at most one point
which is not isolated, and if so it has a neighborhood base
consisting of a decreasing sequence
of (open) sets of length $\theta$. \newline
This suffices to show that $X$ is $(< \lambda)$-metrizable when
$\theta = \aleph_0$ and as required generally (for 4)).

As for showing that $X$ is not $CWH$ (hence not metrizable and not normal),
note that $\{ y_i:i < \lambda \} \cup
\{ z_h:h \in H \}$ is a discrete subspace. \newline
If it is separated, we have a sequence of pairwise disjoint neighborhoods:
\newline
$\langle u_{g(i)}[y_i]:i < \lambda \rangle \char 94 \langle
u_{\zeta(h)}[z_h]:h \in H \rangle$.  But $H$ is not free (in the sense of
Definition 4.1.A)) and we get a contradiction.
\smallskip

\noindent
\underbar{$(b)^+ \Rightarrow (b)$}. \newline
Trivial. \newline
\smallskip

\noindent
\underbar{$(b) \Rightarrow (b)^+$}. \newline
Let $X$ exemplify the second clause so without loss
of generality $|X| = \lambda$.  Let $Y$ be a discrete subspace of cardinality
$\lambda$ which cannot be separated. 
Let $X^+$ be the topology $X$ on the set of points of $X$ generated by 
basic open sets of $X$ and 
$\{ \{ x \}:x \in X \backslash Y \}$.

Now $X^+$ is not $\lambda-CWH$ ($Y$ still exemplifies it).  But $X^+$ is
$(< \lambda)$-metrizable as:

If $Z \subseteq X,|Z| < \lambda$, then we can find a sequence
$\langle u_z:z \in Z \cap Y \rangle$ of pairwise disjoint open sets, and in
$X \restriction u_z$, every point is isolated except $z$, which has a
neighborhood basis of cardinality $\aleph_0$, and every
$x \in Z \backslash \dsize \bigcup_{z \in Z \cap Y} u_z$ is isolated. \newline
This is enough.
\smallskip

\noindent
\underbar{$(b)^+ \Rightarrow (c)^+$}  \newline
Trivial  (as $(< \lambda)$-metrizable $\Rightarrow (< \lambda)-^*CWN$).
\smallskip

\noindent
\underbar{$(c)^+ \Rightarrow (c)$}  \newline
Trivial.
\smallskip

\noindent
\underbar{$(c) \Rightarrow (b)^+$}  \newline
\medskip

If $X,\langle Y_i:i < \alpha \rangle$ exemplifies clause (c) in (1)
with $\langle u_\zeta(y):\zeta < \theta \rangle$ a decreasing 
neighborhood basis of $y$; we can get another example $X'$ to the
third clause, as follows. \newline
We are, without loss of generality, assuming that $|X| = \lambda$.  Then

$$
\align
X' = \dsize \bigcup_{i < \alpha} Y_i \cup \biggl\{
x_{y,z,\zeta,\xi}:\,&\text{for some } i \ne j < \alpha,y \in Y_i,\\
  &z \in Y_j,u_\zeta[y] \cap u_\zeta[z] \ne \emptyset \biggr\}
\endalign
$$

\noindent
with the neighborhood bases for $y,z \in \dsize \bigcup_{i < \alpha} Y_i$
given by

$$
\align
u'_\varepsilon[t] = \{ y \} \cup \biggl\{ x_{y,z,\zeta,\xi}:
&x_{y,z,\zeta,\xi} \in X',\\
  &t = y \wedge \varepsilon \le \zeta \text{ or } t = z \wedge \varepsilon
\le \xi \biggr\}
\endalign
$$

\noindent
and $x_{y,z,\zeta,}$ isolated. \newline
Clearly $Y =: \dsize \bigcup_{i < \alpha} Y_i$ is discrete.  Assume that
$\langle u'_{\varepsilon(y)}[y]:y \in Y \rangle$ is a sequence of pairwise
disjoint open sets.  \underbar{Then} let

$$
U_i = \cup\{ u_{\varepsilon(y)}[y]:y \in Y_i \}.
$$

\noindent
So in $X,U_i$ is an open set (as a union of open sets),

$$
Y_i \subseteq U_i \text{ as } y \in u_{\varepsilon(y)}[y]
$$

$$
\align
i \ne y \Rightarrow U_i \cap U_j = \emptyset &\Rightarrow \exists y \in
U_i,\exists z \in U_j(U_{\varepsilon(y)} \cap U_{\varepsilon(z)} \ne
\emptyset) \\
  &\Rightarrow X_{y,\zeta,\varepsilon(y),\varepsilon(\zeta)}
\text{ is well defined} \\
  &\Rightarrow \text{ in } X' \text{ we have that } U'_{\varepsilon(y)}
\cap U'_{\varepsilon(z)} \ne \emptyset.
\endalign
$$

\noindent
This is a contradiction.
\medskip

So we conclude that $y$ cannot be separated in $X'$, so $X'$ is not
$X-CWH$.

Next, assume that $Z \subseteq X',|Z| < \lambda$, so in $Z,\langle Y_i \cap
Z:i < \alpha,y_i \cap z \ne \emptyset \rangle$ can be separated, say by
$\langle U_i:i < \alpha,Y_i \cap Z \ne \emptyset \rangle$.  So for
$y \in Y \cap Z$, there is an $\varepsilon(y)$, such that
$U'_{\varepsilon(y)}[y] \in U_i$ (the isolated points in $X' \cap Z
\backslash Y$ can be taken care of easily so we ignore them).

Now, if $Y_1 \ne z \in Y \cap Z$ then:

$$
\text{if } (\exists i)(y,z \in Y_i) \text{ then } u'_{\varepsilon(y)} \cap
u'_{\varepsilon(z)} = \emptyset
$$

\noindent
(for any choice of $\varepsilon(y),\varepsilon(z)$ if)

$$
y \in Y_i,z \in Y_j,i \ne j, \text{ if }
u'_{\varepsilon(y)} \cap u'_{\varepsilon(z)} \ne \emptyset
$$

\noindent
then $X_{y,z,\varepsilon(y),\varepsilon(z)}$ exists, so

$$
\emptyset \ne u_{\varepsilon(y)} \cap U_{\varepsilon(z)} \subseteq U_i \cap
U_j
$$

\noindent
which is a contradiction. \newline
That $X'$ is $(< \lambda)$-metrizable now follows as in $(b) \Rightarrow
(b)^+$.
\medskip

\noindent
\underbar{$c) \Rightarrow d)$}. \newline
Assume that $X$ is a normal first countable $(< \lambda)-{}^*CWN$ not
$\lambda-^*CWN$-space, without loss of generality with the set of points
$\lambda$, so there is a
sequence $\langle Y_i:i < \alpha \rangle$ of pairwise disjoint subsets of
$X,Y_i \ne \emptyset$, $Y_i$ is clopen in $X \restriction
(\dsize \bigcup_{j<\alpha} Y_j)$ and $\langle Y_i:i < \alpha \rangle$
cannot be separated.  For $y \in Y =:\dsize \bigcup_{i< \lambda} Y_i$ let
$\bar u Y = \langle u_\zeta[y]:\zeta < \theta \rangle$ be a neighborhood
basis of the topology for $y$,
and without loss of generality $\varepsilon < \zeta < \theta
\Rightarrow u_\zeta[y] \subseteq u_\varepsilon[y]$.  Let

$$
\align
H = \biggl\{ (h,\bar u):&\text{ for some } i < \alpha \text{ and }
\text{ for some } y \in Y_i, \\
  &(k,\bar u) = (y,\bar u_y), \text{ which means}: \\
  &\text{Dom}(h) =
\dsize \bigcup_{j \ne i} Y_j,c \ell(h(z)) = \{ (\zeta,\xi) \in \theta
\times \theta:u_\zeta[y] \cap u_\xi[z_z] = \emptyset\} \\
  &\text{and } \bar u \text{ is } \langle u_\xi(y) \cap \text{ Dom}(h):
\zeta < \theta \rangle \biggr\}.
\endalign
$$

\noindent
Note that $h(z)$ is uniquely determined by $c \ell(h(z)),c \ell_Q(k(z))$.
As we check that $H$ exemplifies $SQw_{\aleph_0}$, i.e. the clauses
in 4.1A(3).
Clauses (a), (b) are immediate.  As for clause (c), let $H' \subseteq H,
|H'| < \lambda$, and \newline
$Z' \subseteq \cup \{ \text{Dom}(h):(h,\bar u) \in H \},
|Z'| < \lambda$, let \newline
$Y' =: \{ y:y \in \dsize \bigcup_{i < \alpha} Y_i$, and
$y \in Z'$ or $(h_y,\bar u_y) \in H' \}$, so $|Y'| < \lambda$; we can find
\newline
$X' \subseteq X,|X'| \le |Y'| + \theta < \lambda$ such that
$Y' \subseteq X'$, and for every $y,z \in Y',\zeta < \theta,\xi < \theta$,
we have
$u_\zeta[y] \cap u_\xi[z] \ne \emptyset \Rightarrow u_\zeta[y] \cap
u_\xi[z] \cap X' \ne \emptyset$.  As $|X'| < \lambda$ we know that $X'$ (i.e.
$X \restriction X'$) is $CWN$, and $\langle Y_i \cap X':i < \alpha \rangle$
is a discrete sequence of closed sets in $X'$ hence there is a function
$g:Y' \rightarrow \theta$ such that
\medskip
\roster
\item "{$(*)$}"  if $i < j < \alpha,
y \in Y' \cap Y_i,z \in Y' \cap Y_j$, then \newline
$u_{g(y)}[y] \cap u_{g(z)}[z] = \emptyset$ (intersecting with $X'$ is
immaterial).
\endroster
\medskip

\noindent
Hence by the choice of $g$
\medskip
\roster
\item "{$(**)$}"  if $i \ne j(i < \alpha,j < \alpha), y \in Y' \cap Y_i,
z \in Y' \cap Y_j$ \newline
then $(g(y),g(z)) \in c \ell(h_y(z))$. 
\endroster
\noindent
This is enough.
\medskip

We are left with proving that $H$ is not free, so suppose
$f,g:Y \rightarrow \theta$ satisfies
\medskip
\roster
\item "{$\bigotimes$}"  for every $y \in Y$, \newline
for every $z \in \text{ Dom}(h_y),(g(y),f(z)) \in c \ell (h_y(z))$,
\endroster
\medskip

\noindent
so without loss of generality $f = g$. \newline
For $i < \alpha$ let

$$
U_i = \cup \{ U_{g(y)}[y]:y \in Y_i \}.
$$

\noindent
So $U_i$, being the union of open sets is open. \newline
If $i < j, y \in Y_i, z \in Y_j$ then

$$
\align
u_{g(y)}[y] \cap u_{g(y)}[z] \ne \emptyset &\Rightarrow (g(y),g(z)) \in 
c \ell(h_y(z)) \\
  &\Rightarrow (g(y),f(z)) = (g(y),g(z)) \in c \ell(h_y(z)).
\endalign
$$

\noindent
Contradiction, by the choice of $f$ and $g$. \newline
So $u_{g(y)}[y] \cap u_{g(z)}[z] = \emptyset$, as $y \in Y_i, z \in Y_j$ were
arbitrary, $U_i \cap U_j = \emptyset$. \newline
We conclude that $\langle Y_i:i < \lambda \rangle$ can be separated, which
is a contradiction. \newline
2)  We prove each implication \newline
(A)  $\lambda \in SQd_{\theta,\sigma} \Rightarrow \lambda \in SQ_
{\theta,\sigma} \Rightarrow \lambda \in SQw_{\theta,\sigma}$.  Obvious.
\newline
(B) $\lambda \in SQw_{\theta,\sigma} \Rightarrow SQd_{\theta,\sigma}
\cap [\lambda,\lambda^\theta] \ne \emptyset$ when $\sigma \le \theta^+$.
\newline

Assume that $H$ exemplifies $\lambda \in SQw_{\theta,\sigma}$.  By the
definition $(h,\bar u) \in H \Rightarrow u_\zeta \subseteq \text{ Dom}(h)
\and \dsize \bigcap_{\xi < \theta} u_\xi = \emptyset$.  Let for
each $(h,\bar u) \in H$,

$$
\align
H^*_{(h,\bar u)} = \biggl\{f:&\,f \text{ is a function from ordinals to }
Pie(\theta \times \theta) \text{ and } \text{ Dom}(f) = \nu \\
  &\text{ for some set } v, v \subseteq \text{ Dom}(h),|v| = \theta, \\
  &\text{ but } \zeta < \theta \Rightarrow |v \backslash u_\zeta| < \theta,
\text{ and } \\
  &(\forall \alpha \in \nu)[c \ell(f(\alpha)) \supseteq c \ell(h(\alpha))],
\text{ and } f \text{ is simple} \biggr\}
\endalign
$$

\noindent
and $H^* = \cup \{H^*_{(h,\bar u)}:(h,\bar u) \in H \}$.
\newline

\noindent
It is easy to check that $H^*$ satisfies clauses (a) and (b) from 4.1A(1)
and (e) of 4.1A(2) and $|H^*| = \lambda^\theta$.
\medskip

As for clause (c) of 4.1A(1), let $H' \subseteq H,|H'| < \lambda$, let
$H' = \{ f_j:j < j(*) \},j(*) < \lambda$, and $(h_j,v_j)$ as
in the definition of $H^*_{(h,\bar u)}$ for some $(h_j,\bar u_j) \in H$.
Define $H'' = \{ (h_j,\bar u_j):j < j(*) \}, Y = \dsize \bigcup
_{j < j(*)} v_j$.  Now $H''$ is a subset of $H$ of cardinality $< \lambda$,
$Y \subseteq Ord$ and $|Y| < \lambda$ so as $H$ exemplifies
$\lambda \in SQw_{\theta,\sigma}$, we can find a $\langle g_i:i < i(*)
\rangle$, $i(*) < \sigma, \, g_i \in {}^\lambda \theta$ and for every
$(h_j,\bar u_j) \in H''$ for some $i = i(j) < i(*)$ we have \newline
$(\exists \zeta < \theta)(\exists
\xi < \theta)(\forall \alpha \in u_{j,\zeta} \cap Y)
[(g_i(\alpha),\xi) \in c \ell(h_j(\alpha))]$. \newline

Now $\langle g_i:i < i(*) \rangle$ are O.K. for $H'$, too, as
$c \ell(f_j(\alpha)) \supseteq c \ell(h_j(\alpha))$ and
$|v_j \backslash u_{i,\zeta}| < \theta$.
\medskip

We are left with clause (d) of 4.1A(1), so assume $i(*) < \sigma$
and $g_i \in
{}^\lambda \theta$ for $i < i(*)$ exemplifies $H^*$ is $\sigma$-free.
By the choice of $H$ for some \newline
$(h,\bar u) \in H$ we have
$\dsize \bigwedge_i \neg(\exists \zeta < \theta)(\exists \xi
< \theta)(\forall \alpha \in u_\zeta)[(g_i(\alpha),\xi) \in c
\ell(h(\alpha))]$. \newline

Let $\langle a_i:i < i(*) \rangle$ be a partition of $\theta$ to unbounded
subsets, and we choose by induction on $\zeta < \theta$, an ordinal
$\alpha_\zeta \in u_\zeta$ and $\gamma_\xi < \theta$ such that if
$\alpha_\zeta \in a_i$ then

$$
\Upsilon_\zeta \in \theta \backslash \left[
\dsize \bigcup \Sb \xi < \zeta \\ i < \sigma \endSb
(g_i(\Upsilon_\xi) \cup \Upsilon_\varepsilon) + 1 \right]
$$

$$
(g_i(\alpha_\zeta),\Upsilon_\zeta) \notin c \ell(h(\alpha_\zeta))
$$

\noindent
and let $f(\alpha_\zeta)$ be such that

$$
c \ell(f(\alpha_\zeta)) = \{(\gamma_1,\gamma_2):\gamma_1 < \theta,
\gamma_2 < \theta, \text{ and }
(\gamma_1,\gamma_2) \nleq (g_i(\alpha_\zeta),\Upsilon_\zeta) \}.
$$

\noindent
Let $v =: \{ \alpha_g:\zeta < \theta \}$, so $f \in H^*(h,\bar u) \subseteq
H^*$ exemplifies that $\langle g_i:i < i(*) \rangle$ exemplify tht $H^*$ is
$\sigma$-free.  We can finish by 4.3A(1). \newline
(3)  As in 2), $\lambda \in SPd_{\theta,\sigma} \Rightarrow
\lambda \in SP_{\theta,\sigma} \Rightarrow \lambda \in SPw_{\theta,\sigma}$
is obvious. \newline
We need to prove that $\lambda \in SPw_{\theta,\sigma} \Rightarrow
SPd_{\theta,\sigma} \cap[\lambda,\lambda^\theta] \ne \emptyset$ when
$\sigma \le \theta^+$. \newline
The proof if similar to that of (2).  We start with $H$ exemplifying that
$\lambda \in SPw_{\theta,\sigma}$.  We assume that for each
$(h,\bar u) \in H,\bar u$ is standard.  So for $(h,\bar u) \in H$, we
define

$$
\align
H^*_{(h,\bar u)} = \biggl\{
f:\,&f \text{ is a function from ordinals to } \theta \text{ and }
f \text{ is } 1-1, \\
  &\text{ and for some set } v \subseteq Dom(h), \text{ we have that }
|v| = \theta, \text{ but } \\
  &\zeta < \theta \Rightarrow |v \backslash u_\zeta| < \theta, \text{ while }
(\forall \alpha \in v)f(\alpha) \le h(\alpha) \biggr\}.
\endalign
$$

\noindent
Let $H^* = \cup \{ H^*_{(h,\bar u)}:(h,\bar u) \in H \}$. \newline
Checking that this $H^*$ is as required is similar to (2).  For example, to
see 4.1.1)d), suppose that $i(*) < \sigma$ and $\langle g_i:i < i(*) \rangle$
exemplify that $H^*$ is free.  By the choice of $H^*$, there is an
$(h,\bar u) \in H$ such that

$$
\dsize \bigwedge_i \neg(\exists \zeta < \theta)(\exists \zeta < \theta)
(\forall \alpha \in u_\zeta)[h(\alpha) \le \text{ max}\{ g_i(\alpha),\xi \}].
$$

\noindent
Let $\langle a_i:i < i(*) \rangle$ be as in (2), and we choose by induction
on $\zeta < \theta$, an ordinal $\alpha_\zeta \in u_\zeta$ and $\Upsilon
_\zeta \in \theta$ such that

$$
\alpha_\zeta \in a_i \Rightarrow \Upsilon_\zeta \in \theta \backslash
\left[ \dsize \bigcup \Sb i < \sigma \\ \xi < \zeta \endSb
(g_i(\gamma_\xi) \cup \gamma_\xi) + 1 \right]
$$

\noindent
and

$$
h(\alpha_\zeta) > \text{ max}\{g_i(\alpha_\zeta),\gamma_\zeta \}.
$$

\noindent
Then we let $f(\alpha_\zeta)$ be such that

$$
f(\alpha_\zeta) \le \text{ max}\{g_i(\alpha_\zeta,\gamma_\zeta \}
$$

\noindent
but

$$
f(\alpha_\zeta) \notin \{f(\alpha_\xi):\xi < \zeta \}.
$$

\noindent
4) Included in the proof of (1). \hfill$\square_{4.4}$
\enddemo
\bigskip

\proclaim{4.5 Claim}  Assume $\lambda \in SP_{\theta,\sigma}$,
 $\mu$ is a strong limit with
$cf(\mu) > \theta$, and \newline
$2^\mu = \mu^+ > \lambda$.

Then there is a $\kappa \in [\lambda,\mu^+]$, a regular cardinal such that
$\kappa \in SP^+_{\theta,\sigma}$ where
\endproclaim
\medskip

\definition{4.6 Definition}  1)  $\kappa \in SP^+_{\theta,\sigma}$ means that
$\kappa$ is regular $> \theta$ and we can
find an $S \subseteq \{ \delta < \kappa:cf(\delta) = \theta \}$ stationary,
$\bar \eta = \langle \eta_\delta:\delta \in S \rangle$, $\bar h = \langle
h_\delta:\delta \in S \rangle$, such that
\medskip
\roster
\item "{(a)}"  $\eta_\delta$ is a strictly increasing sequence
 of ordinals \newline
 of length $\theta$ with limit $\delta$
\item "{(b)}"  $h_\delta:\text{Rang}(\eta_\delta) \rightarrow \theta$ is
strictly increasing
\item "{(c)}"  $H = \{ h_\delta:\delta \in S \}$ is ($< \kappa$)-$\sigma$-
free not $\sigma$-free (in 4.1's sense).
\endroster
\medskip

\noindent
2)  $\kappa \in SP^*_{\theta,\sigma}$ if in the above we add:
\roster
\item "{(d)}"  $h_\delta(\eta_\delta(\varepsilon))$ depend on
$\eta_\delta(\varepsilon)$ only
\item "{(e)}"  $\bar \eta$ is tree like, i.e. $\eta_{\delta_1}
(\varepsilon_1) = 
\eta_{\delta_2}(\varepsilon_2) \Rightarrow \varepsilon_1 = \varepsilon_2
 \and \eta_{\delta_1}\restriction \varepsilon_1 = \eta_{\delta_2}
\restriction \varepsilon_2$.
\endroster
\enddefinition
\bigskip

\remark{Remark}  The assumption $``\mu^+ = 2^\mu"$ (in 4.5) is very reasonable
because of 4.2(2) (and 4.2(3) from the topological point of view).
\endremark
\bigskip

\demo{4.6A Observation}  1) $SP^*_{\theta,\sigma}
 \subseteq SP^+_{\theta,\sigma} \subseteq SP_{\theta,\sigma}$. \newline
2)  If $\langle h_\delta,\eta_\delta:\delta \in S \rangle$, $\kappa$
satisfies the preliminary requirements and clauses (a), (b) of 4.5 and $H$
is $(< \kappa_1)$-free, $\kappa_1 > \theta$ then for some $\mu \in [\kappa_1,
\kappa]$, $\mu \in SP^+_\theta$. \newline
3)  Similarly for $SP^*_{\theta,\sigma}$.
\enddemo
\bigskip

\demo{Proof} Like 2.3 or 2.4. \hfill$\square_{4.6A}$
\enddemo
\bigskip

\demo{4.6.B Conclusion}  For $\lambda > \theta = cf(\theta),\chi =
\beth_\chi > \lambda$, the following
are equivalent:
\roster
\item "{(a)}"  for some $\mu \in [\lambda,\chi)$, $\mu \in SP_\theta$
\item "{(b)}"  for some $\mu \in [\lambda,\chi)$, $\mu \in SP^+_\theta$.
\item "{(c)}"  for some $\mu \in [\lambda,\chi),\mu \in SP^*_\theta$.
\endroster
\enddemo
\bigskip

\demo{Proof}  By 4.5, (b) $\Rightarrow$ (a), as for (a) $\Rightarrow$ (b),
let $\mu = \beth_{\lambda^{+(\theta^+)}}$, if $pp(\mu) > \mu^+$ use 4.2(2) and
if $pp(\mu) = \mu^+$ use 4.4. \hfill$\square_{4.6.B}$
\enddemo
\bigskip

\demo{Proof of 4.5}  Use $\diamondsuit_{\{ \delta < \mu^+:cf(\delta) = \theta
\}}$ and imitate 3.3. \hfill$\square_{4.5}$
\enddemo
\bigskip

\proclaim{4.7 Claim}  Assume $\theta = \theta^{< \theta}$ or
$\exists F \subseteq {}^\theta \theta$ which is cofinal in
${}^\theta \theta$ and $| \{ f \restriction \zeta:
f \in F,\zeta < \theta \}| \le \theta$.
Let $\langle h_\delta,\eta_\delta:\delta \in S
\rangle$ exemplify $\lambda \in SP^*_\theta$ (even omitting ``$\eta_\delta$
converge to $\delta,\eta_\delta$ strictly increasing").
Then any $\theta^+$-complete
forcing preserves the non-freeness of \newline
$\{ h_\delta:\delta \in S \}$.
\endproclaim
\bigskip

\demo{Proof}  Instead of the domain of the functions $h_\delta$
being a subset of $\lambda$, we can assume that it is 
$T = \{ \eta_\delta \restriction \zeta:\delta \in S,\zeta < \theta
\text{ a successor}\}$ (identify $\eta_\delta(\zeta)$ with $\eta_\delta
\restriction(\zeta + 1)$, so $\text{ Dom}(h_\delta) = \{ \eta_\delta
\restriction \zeta:\zeta < \theta \text{ is a successor ordinal } \}$).
Suppose $Q$ is
a $\theta^+$-complete forcing notion, $p \in Q$ and $p \Vdash ``{\underset\sim
\to g}:T \rightarrow \theta$ exemplifies $\{ h_\delta:\delta \in S\}$ is free
".  We now define by induction on $\ell g(\eta) < \theta$ a seuqnce
$\langle p_{\eta,t},\varepsilon_{\eta,t},v:t \in T_\eta \rangle$ for
$\eta \in T$ such that:
\medskip
\roster
\item "{$(\alpha)$}"  $T_\eta \subseteq {}^{\ell g(\eta) \ge}\theta$, is
closed under initial segments
\item "{$(\beta)$}"  $t \triangleleft s \in T_\eta \Rightarrow p_{\eta,t}
\le_Q p_{\eta,s}$
\item "{$(\gamma)$}"  if $t \in T_\eta$, either $\dsize \bigwedge_
{\zeta < \theta} t \char 94 <\zeta> \in T_\eta$ or $\dsize \bigwedge_
{\zeta < \theta} t \char 94 <\zeta> \notin T_\eta$
\item "{$(\delta)$}"  If $t \in {}^{\ell g(\eta) \ge} \theta$, $\ell g(t)$
is a limit ordinal and $(\forall \zeta < \ell g(t))(t \restriction \zeta \in
T_\eta)$, then $t \in T_\eta$,
\item "{$(\varepsilon)$}" if $\nu \triangleleft \eta$ then $T_\nu \subseteq
T_\eta$ and $t \in T_\nu \Rightarrow (p_{\eta,t},\varepsilon_{\eta,t}) =
(p_{\nu,t},\varepsilon_{\nu,t})$
\item "{$(\zeta)$}"  if $\ell g(\eta)$ is a limit ordinal then \newline
$T_\eta = \{ t:t \in \dsize \bigcup_{\nu \triangleleft \eta} T_\nu$ or
$\ell g(t)$ is a limit ordinal and $(\forall s)[s \triangleleft t
\Rightarrow s \in \dsize \bigcup_{\nu \triangleleft \eta} T_\nu] \}$.
\item "{$(\eta)$}"  assume $\eta = \nu \char 94 <\alpha>$ and $s$ is a
 $\triangleleft$-maximal element of $T_\nu$, then:
{\roster
\itemitem{ (a) }  if $\{ \zeta < \theta:p_{\eta,s} \nvDash_Q
 ``{\underset\sim\to g}(\eta) \ne \zeta" \}$ is bounded in $\theta$ \newline
then $s$ is a $\triangleleft$-maximal element of $T_\eta$.
\itemitem{ (b) }  if $A = \{ \zeta < \theta:p_{\eta,s} \nvDash
``{\underset\sim\to g}(\eta) \ne \zeta" \}$ is unbounded in $\theta$, then for
every $\zeta < \theta$, $s \char 94 <\zeta>$ is a maximal member of $T_\eta$,
and $p_{s \char 94 <\zeta>}$ forces a value $\varepsilon_{s \char 94 <\zeta>}
 > \ell g(\eta)$ to ${\underset\sim\to g}(\eta)$.
\endroster}
\endroster
\medskip

\noindent
We can carry this definition.
\roster
\item "{$(*)$}"  if $\delta \in S$ then for some $\zeta = \zeta_\delta <
\theta$ and $t = t_\delta \in T_{\eta_\delta \restriction \zeta}$ we have:
$t$ is a $\triangleleft$-maximal member of $T_{\eta \restriction \xi}$ for
every $\xi \in [\zeta,\theta)$.
\endroster
\medskip

\noindent
[why?  otherwise we can construct a $t \in {}^\theta \theta$ such
that $(\forall s)[s \triangleleft t \Rightarrow s \in \dsize \bigcup
_{\xi < \theta} T_{\eta_\delta \restriction \xi} ]$, $t(\varepsilon) >
\varepsilon$
and for unboundedly many $\xi < \theta$, for some $s \char 94 <\zeta>
\triangleleft \, t$ we have \newline
$s \char 94 <\zeta> \in T_{\eta_\delta
\restriction (\xi+1)} \backslash T_{\eta_\delta \restriction \xi}$,
$\varepsilon_{s \char 94 <\zeta>} >
h_\delta(\eta_\delta \restriction (\xi + 1)),\xi$. \newline
Now $\{ p_{\nu,s}:\nu \triangleleft \eta_\delta,s \triangleleft t, s \in
T_\nu \}$ has an upper bound in $Q$, $p^*$.  Then $p^*$ forces for
$\underset\sim\to g (\eta \restriction (\xi + 1))$ a value $> h_\delta
(\eta \restriction (\xi + 1)),\xi$; this is a contradiction to \newline
$p \Vdash
``\underset\sim\to g$ exemplifies the freeness of $\{ h_\delta:
\delta \in S \}$"]. \hfill$\square_{4.7}$
\enddemo
\bigskip

\proclaim{4.8 Theorem}  Assume $\lambda < \mu,(\forall \kappa < \mu)[\chi
^{\aleph_0} < \mu)$ (possibly $\mu = \infty$).  Then the following are
equivalent:
\medskip
\roster
\item "{(A)}"  There is a space $X$ such that:
{\roster
\itemitem{ (a) }  $X$ is $(< \lambda)$-metrizable
\itemitem{ (b) }  $X$ is not metrizable
\itemitem{ (c) }  $X$ has $< \mu$ points.
\endroster}
\item "{(B)}"  There is a first countable Hausdorff space $X$ such that:
{\roster
\itemitem{ (a) }  $X$ is $(< \lambda)$-CWH
\itemitem{ (b) }  $X$ is not $\lambda$-CWH
\itemitem{ (c) }  $X$ has $< \mu$ points.
\endroster}
\item "{(B)$^+$}"  There is a space $X$ like in (B), and in addition
{\roster
\itemitem{ $(a)^+$ }  $X$ is $(< \lambda)$-metrizable.
\endroster}
\item "{(C)}"  There is a first countable Hausdorff space $X$ such that:
{\roster
\itemitem{ (a) }  $X$ is $(< \lambda)-^*$ CWN
\itemitem{ (b) }  $X$ is not $\lambda-^*$ CWN
\itemitem{ (c) }  $X$ has $< \mu$ points.
\endroster}
\item "{(C)$^+$}"  There is an $X$ like in (C), and in addition,
{\roster
\itemitem{ (a)$^+$ }  $X$ is $(< \lambda)$-metrizable.
\endroster}
\item "{(D)}"  there is a family $H$ of functions with domains countable
sets of ordinals and range $\subseteq \omega$ such that:
{\roster
\itemitem{ (a) }  $H$ is $(< \lambda)$-free
\itemitem{ (b) }  $H$ is not free
\itemitem{ (c) }  $|H| < \mu$.
\endroster}
\item "{$(D)^+$}"  as in (D) and
{\roster
\itemitem{ (d) }  $\cup \{ \text{Dom}(h):h \in H \} = \lambda' \in [\lambda,
\mu)$
\itemitem{ (e) }  each $h$ is one to one.
\endroster}
\item "{$(D)'$}"  $[\mu,\lambda) \cap SP_{\aleph_0} \ne \emptyset$
\item "{$(D)''$}"  $[\mu,\lambda) \cap SPw_{\aleph_0} \ne \emptyset$
\item "{$(D)'''$}"  $[\mu,\lambda) \cap SPd_{\aleph_0} \ne \emptyset$
\item "{(E)}"  there is a ${\overset =\to u} = \langle <u_{\alpha,n}:
n < \omega >:\alpha \in v \rangle$,$u_{\alpha,n+1} \subseteq
u_{\alpha,n} \subseteq v$, such that:
{\roster
\itemitem{ (a) }  ${\overset =\to u}$ is not free
\itemitem{ (b) }  for $v' \in [v]^{< \lambda}$, ${\overset =\to u}
\restriction v'$ is free
\itemitem{ (c) }  $|v| < \mu$.
\endroster}
\item "{$(E)'$}"  $[\mu,\lambda) \cap SQ_{\aleph_0} \ne \emptyset$
\item "{$(E)''$}"  $[\mu,\lambda) \cap SQw_{\aleph_0} \ne \emptyset$
\item "{$(E)'''$}"  $[\mu,\lambda) \cap SQd_{\aleph_0} \ne \emptyset$
\endroster
\endproclaim
\bigskip

\proclaim{4.8A Theorem}  In 4.8 if $(\forall \kappa < \mu)(\beth_{\theta^+}
(\kappa) < mu)$ (really $(\forall \kappa < \mu)(\beth_{\omega_1}(\kappa)
< \mu)$ is O.K.  Equivalently $\mu = \beth_\delta = \theta^+(\delta)$ then
we can add
\medskip
\roster
\item "{$(F)$}"  for some regular $\kappa \in [\lambda,\mu)$ we have
$INCWH(\kappa)$
\item "{$(F)'$}"  $\lambda \in SP^+_{\aleph_0}$
\item "{$(F)''$}" $\lambda \in SP^*_{\aleph_0}$
\endroster
\endproclaim
\bigskip

\demo{Proof of 4.8A}  By 4.4(1) (the $(b) \Leftrightarrow (b)^+
\Leftrightarrow (c) \Leftrightarrow (c)^+$ part) we know the equivalence
of (A), (B), $(B)^+$, (C) $(C)^+$. \newline
By x.x $(D) \Leftrightarrow (D)'$. \newline
By 4.4(3) we have $(D)' \Rightarrow (D)'' \Rightarrow (D)'''$. \newline
By x.x $(E) \Leftrightarrow (E)'$. \newline
By 4.4(3A) $(E)' \Rightarrow (E)'' \Rightarrow (E)''$. \newline
By 4.1A(2) $(E)'' \Rightarrow (D)'$. \newline
By 4.1A(1) $(D)' \Rightarrow (E)',(D)'' \Rightarrow (E)'',
(D)''' \Rightarrow (E)''$. \newline
\medskip

Together we get the equivalence of $(D), (E), (D)', (E)', (D)'', (E)'',
(D)'', (E)''$.

By 4.4(1) $(E)' \Rightarrow (A) \Rightarrow (E)'''$, so by the last
sentence and the first paragraph we have finished the proof of 4.8.  For
4.8A use 4.6B.
\enddemo
\bigskip

\demo{4.4 Fact}  Let $\lambda = cf(\lambda) > \theta = cf(\theta)$.
\medskip
\roster
\item "{(A)}"  There is $H = \dsize \bigcup_{i < \lambda} H_i$ such that:
{\roster
\itemitem{ $(\alpha)$ }  $H_i$ is increasing continuous
\itemitem{ $(\beta)$ }   $H_\theta$ is a family of functions $h$,
$Dom(h)$ is a set of $\theta$ ordinals, $h$ is one to one
\itemitem{ $(\gamma)$ }  each $H_i$ is free, but $H$ is not free.
\endroster}
\endroster
\medskip

\noindent
$(B) = (B)_{r,\theta}$.  Let  $X = X_{\lambda,\theta} =: {}^\lambda \theta$

$$
\align
F = F_{\lambda,\theta} =: \biggl\{
f: \, &f \text{ a partial function from } X \text{ to } \theta,|Dom \, f|
= \theta \text{ and} \\
  &(\forall^*i < \lambda)(\forall^* \eta \in \text{ Dom}(f))
[f(\eta) \le \eta(i)] \\
  &\text{ and } f \text{ is one to one } \biggr\}.
\endalign
$$

\noindent
Then there is no $G:X \rightarrow \omega$ such that

$$
f \in F \Rightarrow f(\forall^*\eta \in \text{ Dom}(f))
[f(\eta) \le G(\eta)].
$$

\noindent
Then \newline
$(A) \Leftrightarrow (B)$.  
\enddemo
\bigskip

\demo{Proof} \underbar{$(A) \Rightarrow (B)$}. \newline

Let $H,H_i(i < \lambda)$ exemplifies (A), let $A = \cup \{ Dom(h):h \in H \}$,
and let $g_i:A = \theta$ exemplify ``$H_i$ is free". \newline
We define an equivalence relation $E$ on $A:\alpha \exists_\beta
\Leftrightarrow \dsize \bigwedge_{i < \lambda} g_i(\alpha) = g_i(\beta)$.
If for some $h \in H$ and $\alpha,(\alpha/E) \cap \text{ Dom}(h)$ has
cardinality $\theta$, choose $i < \lambda$ such that $h \in H_i$, and
$g_i$ cannot satisfy the requirement.  let $h^\otimes$ be a function with
domain $\text{Dom}(h)$, $h^\otimes(\alpha) = \text{ sup}\{ h(\beta):\beta
\in \alpha/E \}$.  Now $H' =: \{ h^\otimes:h \in H\},H''_i \{ h^\otimes:
h \in H'_i \}$ exemplifies (A) too.  So without loss of generality $E$ is
the equality on $A$. \newline
Next for each $\alpha \in A$ let $\eta_\alpha \in {}^\lambda \theta(=X)$ be
defined by $\eta_\alpha(i) = g_i(\alpha)$, so $\alpha \ne \beta \Rightarrow
\eta_\alpha \ne \eta_\beta$.  For $h \in H$ let $\text{Dom}(h) = \{
\alpha_{h,\zeta}:\zeta < \theta \}$ such that $\langle h(\alpha^*_{h,\zeta}):
\zeta < \theta \rangle$ is strictly increasing.  For $h \in H$ let the
function $f_h$ be defined by:

$$
Dom(f_h) = \{ \eta_{\alpha_{h,\zeta}}:\zeta < \theta \},
f_h(\eta_{\alpha_{h,\zeta}}) = h(\alpha_{h,\zeta}).
$$ 

\noindent
Now
\roster
\item "{$(*)$}"  $h \in H \Rightarrow f_h \in F$.
\endroster

\noindent
[Why?  Let $i(*) = \text{ min}\{i:h \in H_i \}$ (well defined as
$H = \dsize \bigcup_{i < \lambda} H_i)$, so $i \in [i(*),\lambda)$ implies
$h \le^* (g_i \restriction \text{ Dom } h)$.  So for some $\zeta(*) <
\theta$, for every $\zeta \in [\zeta(*),\theta)$ we have
$h(\alpha_{h,\zeta}) \le g_i(\alpha_{h,\zeta})$, but
$f_h(\eta_{\alpha_{h,\zeta}}) = h(\alpha_{h,\zeta})$ and $g_i(\alpha_{h,\zeta}
) = \eta_\alpha(i)$ so: for every $i < \lambda$ large enough for all but
$< \theta$ members $\eta = \eta(\alpha_{h,\zeta})$ of $\text{Dom }f_h,
f_h(\eta) = h(\alpha_{h,\zeta}) < g_i(\alpha_{h,\zeta}) = \eta_{\alpha_{h,
\zeta}}(i) = \eta(i)$ as required].

So assume $G$ is a function from $X$ to $\omega$ such that
\medskip
\roster
\item "{$(**)$}"  $f \in F \Rightarrow (\forall^* \eta \in \text{ Dom}(f))
[f(\eta) \le G(\eta)]$
\endroster
\medskip

\noindent
and we should get a contradiction.  let us define $g \in {}^A \theta$ by
$g(\alpha) = G(\eta_\alpha)$.  So for $h \in H$, we have $f_h \in F$ hence
by $(*) + (**)$ for some $\zeta(*) < \theta,\zeta \in [\zeta(*),\theta)
\Rightarrow f_h(\eta_{\alpha_{h,\zeta}}) \le G(\eta_{\alpha_{h,\zeta}})$.
But $f_h(\eta_{\alpha_{h,\zeta}}) = h(\alpha_{h,\zeta})$, and
$g(\alpha_{h,\zeta}) = G(\eta_{\alpha_{h,\zeta}})$ so $\zeta \in [\zeta(*),
\theta) \Rightarrow h(\alpha_{h,\zeta}) \le g(\alpha_{h,\zeta})$.  So
$g$ shows that $H$ is free, contradiction.  We have proved (B).
\enddemo
\bigskip

\noindent
\underbar{$(B) \Rightarrow (A)$}

The demand $A = \dsize \bigcup_{h \in H} \text{ Dom}(h) \subseteq \text{ Ord}$
is immaterial, so let $A = X, H = F_{\lambda,\theta}$.  Lastly for
$i < \lambda$ let $g_i:A \rightarrow \theta$ be $g_i(\eta) = \eta(i)$, and

$$
H_i = \{ f \in F:\text{ for every } j \in [i,\lambda) \text{ we have }
(\forall^* \eta \in \text{ Dom}(f))[f(\eta) \le \eta(i) \}.
$$

\demo{4.10 Conclusion}  $INCWH(\lambda)$ implies
$(B)_{\lambda,\theta}$ of Fact 4.9 implies
$(\exists \mu)[\lambda \le \mu \le 2^\lambda \and INCWH(\mu)]$.
\enddemo
\bigskip

\remark{4.11 Remark}  It is well known that
\medskip
\roster
\item "{$(*)$}"  if there is a real valued measure $m$ on $P(\lambda),
\theta = \aleph_0$ \newline
$G(f) = \text{ Min}\{n:m(f^{-1}(\{ n \}) > 0 \}$
\endroster
\medskip

\noindent
then $G$ contradicts $(B)_{\lambda,\aleph_0}$.
\medskip

Also, it is consistent that $SP_{\aleph_0} \subseteq (2^{\aleph_0})^+$.
This follows from the consistency of the PMEA (Product Measure Extension
Axiom) and Fact 4.9. \newline
The consistency of PMEA is due to Kunen.  See [Fl] for an exposition.
\endremark

\newpage

\bye